\documentclass[12pt]{article}

\usepackage{amssymb}
\usepackage{amsmath}
\usepackage{amsbsy}
\usepackage{amscd}
\usepackage{amsfonts}
\usepackage{amsthm}
\usepackage{mathrsfs}
\usepackage{verbatim}
\usepackage{hyperref}
\usepackage{fullpage}
\usepackage{mathdots}
\usepackage{graphicx,subfigure}
\usepackage[english]{babel}
\usepackage[utf8]{inputenc}
\usepackage{tikz}
\usepackage{adjustbox}
\usepackage{authblk}
\usepackage{caption}

\usepackage{mathtext}
\usepackage{stackengine,scalerel}

\usepackage{tikz-cd}
\usetikzlibrary{backgrounds,fit, matrix}
\usetikzlibrary{positioning}
\usetikzlibrary{calc,through,chains}
\usetikzlibrary{arrows,shapes,automata, petri}

\usepackage[boxsize=1.25em, centerboxes]{ytableau}

\newtheorem{Theorem}[equation]{Theorem}
\newtheorem{Corollary}[equation]{Corollary}
\newtheorem{Lemma}[equation]{Lemma}
\newtheorem{Proposition}[equation]{Proposition}

\theoremstyle{definition}

\newtheorem{Example}[equation]{Example}

\newtheorem{Remark}[equation]{Remark}

\numberwithin{equation}{section}
\numberwithin{figure}{section}

\newcommand{\F}{{\mathbb F}}
\newcommand{\PP}{{\mathbb P}}
\newcommand{\C}{{\mathbb C}}
\newcommand{\Z}{{\mathbb Z}}
\newcommand{\Q}{{\mathbb Q}}
\newcommand{\R}{{\mathbb R}}

\newcommand{\N}{{\mathbb N}}

\newcommand{\mc}[1]{\mathcal{#1}}

\newcommand{\mbf}[1]{\mathbf{#1}}

\newcommand\red[1]{{\color{red}#1}}

\DeclareMathOperator{\Aut}{Aut}

\begin{document}
	
\title{Toric Richardson Varieties}

\author[1,2]{Mahir Bilen Can}
%\author[1]{Mahir Bilen Can}
\author[3]{Pinakinath Saha}

\affil[1,2]{\small{
Theoretical Sciences Visiting Program\\
Okinawa Institute of Science and Technology Graduate University, Onna, 904-0495, Japan
\\

Tulane University, New Orleans, USA, mahirbilencan@gmail.com}}

\affil[3]{\small{Indian Institute of Science, Bangaluru, India, pinakinaths@iisc.ac.in}}

\maketitle

\begin{abstract}
In this article, we provide characterizations of toric Richardson varieties across all types through three distinct approaches: 1) poset theory, 2) root theory, and 3) geometry.
\vspace{.5cm}



\end{abstract}
	
\section{Introduction}

Richardson varieties form an important class of algebraic varieties that arise in many areas of mathematics, including representation theory, intersection theory, and geometric combinatorics. 
As they are formed by intersecting Schubert varieties with opposite Schubert varieties, they are crucial for understanding the intersection theory of flag varieties.
In the present article, instead of using them for intersection theory, we investigate the complexities of the torus actions on Richardson varieties. 
Our primary goal is to characterize, in a type free manner, the Richardson varieties that have complexity zero under the natural action of a maximal torus.

Many fundamental properties of Richardson varieties had been explored by various authors even before Richardson's seminal work~\cite{Richardson1992}. 
Interestingly, the origins of Richardson varieties in Grassmannians can be traced back to Hodge's research in~\cite{Hodge1942}. 
Another significant contribution, more relevant to our paper, comes from Deodhar's work~\cite{Deodhar1985}, where he proved the irreducibility of Richardson varieties in (partial) flag varieties. Deodhar uncovered intriguing connections between the stratifications of Richardson varieties, utilizing them to investigate properties related to Kazhdan-Lusztig theory.

In the early 2000s, Brion and Lakshmibai~\cite{BrionLakshmibai2003} further advanced the `Standard Monomial Theory' using Richardson varieties. 
Concurrently, Marsh and Rietsch~\cite{MarshRietsch2004} applied them to the study of the totally nonnegative part of a flag variety. More recently, Richardson varieties took on an additional combinatorial dimension in the work of Tsukerman and Williams~\cite{TsukermanWilliams2015}, who introduced Bruhat interval posets to explore the geometry of these varieties, albeit primarily in type A. Tsukerman and Williams's work was pushed further by Lee, Masuda, and Park in~\cite{LMP2021}, who characterized toric Richardson varieties (in type A only) in terms of Bruhat interval polytopes. 
Finally, this year, Lee, Masuda, and Park unveiled a wealth of intriguing properties concerning toric Richardson varieties, establishing fascinating connections with Catalan combinatorics in~\cite{LeeMasudaPark2023}. 
Throughout these works, a recurring question that has remained unanswered is the characterization of toric Richardson varieties across all types. Addressing this question served as a primary motivation for our article.

In this context, the central result of our article is as follows.

\begin{Theorem}\label{intro:T1}
Let $G$ be a connected, simple, simply connected algebraic group defined over an algebraically closed field. 
Let $B$ be a Borel subgroup of $G$ containing a maximal torus $T$.
Let $X_w^v$ be a Richardson variety in flag variety $G/B$ such that $\ell(w) - \ell(v) \leq \dim T$. 
Then the following assertions are equivalent:
\begin{enumerate}
\item[(i)] $X_w^v$ is a toric variety with respect to the action of $T$.
\item[(ii)] The Bruhat-Chevalley order on the interval $[v,w]$ does not contain any interval that is isomorphic to symmetric group $\mathbf{S}_3$ with respect to Bruhat-Chevalley order.
\item[(iii)] The interval $[v,w]$ is a lattice. 
\end{enumerate}
\end{Theorem}

It is worth noting that our proofs heavily relied on an influential work by Brenti~\cite{Brenti1994}. 
In this article, Brenti characterizes the intervals $[v,w]$ in condition (ii) by leveraging the $R$-polynomials introduced by Kazdhan and Lusztig. 
In fact, we establish our Theorem~\ref{intro:T1} by combining Brenti's results with another characterization of toric Schubert varieties based on root systems. 
While the root theoretic characterization we discovered is not overly complex to articulate, it does involve a bit more notation than what is presented in Theorem~\ref{intro:T1}. As a result, we have chosen to defer its full exposition to the body of the article.

Theorem~\ref{intro:T1} raises the question of a combinatorial characterization of the singular toric Richardson varieties. 
By extending the result~\cite[Theorem 1.3]{LMP2021} of Lee, Park, and Masuda to simply laced cases, we answer this question for the smooth ones. 
\begin{Theorem}\label{T:Smooth}
Let $G$ be a connected, simple, simply connected algebraic group defined over an algebraically closed field. 
Additionally, we assume that $G$ is of type $A,D$, or $E$. 
Let $X_w^v$ be a toric Richardson variety with respect to the action of a maximal torus $T\subset G$. 
Then the interval $[v,w]$ is a Boolean lattice if and only if $X_w^v$ is a smooth toric variety.
\end{Theorem}

In a concise yet crucial article, Karuppuchammy~\cite{Karuppuchamy2013} proved that a Schubert variety $X_w$ in $G/B$ can be considered a toric variety under the action of $T$ if and only if $w$ is a \red{\em Coxeter type element}, which means it is a product of distinct simple reflections. 
As an application of our two theorems, we are able to establish the following analogous result for toric Richardson varieties in certain Grassmannians.

\begin{Theorem}\label{intro:T2}
Let $G$ be a connected, simple, simply connected algebraic group defined over an algebraically closed field. 
Let $G/P$ be a Grassmannian, where either the parabolic subgroup $P$ is minuscule, or $G$ is of type $G_2$.
Let $T$ be a maximal torus of $P$. 
Then a Richardson variety $X_{wP}^{vP}$ in $G/P$ is a toric variety for the action of $T$ if and only if $w=cv$ for some Coxeter type element $c$.
\end{Theorem}

We are now ready to give a concise overview of our article. 
In the next section, which serves as a preliminary introduction, we establish our notation and provide a brief review of certain tools attributed to Deodhar~\cite{Deodhar1985},  Marsh and Rietsch~\cite{MarshRietsch2004}.
In Section~\ref{S:Characterizations}, we present our theorems, offering various characterizations of toric Richardson varieties in a type-independent manner. In Section~\ref{S:Grassmanns}, we prove Theorem~\ref{intro:T2}.
Finally, in Section~\ref{S:Concluding}, we prove Theorem~\ref{T:Smooth}.

\section{Notation and preliminaries}\label{S:Preliminaries}

In this section, we will setup our notation. 
We follow~\cite{Humphreys} for algebraic groups.

Throughout our article, unless otherwise states, $G$ denotes a connected, simple, simply connected algebraic group defined over an algebraically closed field $k$ of arbitrary characteristic. 
We will use the standard notation, $B,T,W$, and $B^-$ to denote a Borel subgroup, a maximal torus contained in $B$, the Weyl group determined by $G$ and $T$,
and the Borel subgroup opposite to $B$. Note that $B^-$ is the unique Borel subgroup of $G$ such that $B\cap B^- = T$. 
Note that $W$ is defined as the quotient group $N_G(T)/T$, where $N_G(T)$ represents the normalizer of $T$ in $G$. 
For each $w\in W$, we use $n_w$ to represent a chosen representative from $N_G(T)$. 
Occasionally, for ease of notation, we will refer to this representative as $w$.
Note also that $B^{-}\,=\,n_{0}B^{+}n_{0}^{-1}$, where $n_{0}$ is a representative in $N_{G}(T)$ of the longest element $w_{0}$ of $W$. 
For the natural left action of $T$ on $G/B$, the $T$-fixed points are precisely the cosets $e_{w}:= n_wB/B$, where $w\in W$.
For $w\,\in\, W$, the corresponding \red{\em Schubert cell}, denoted $C_{wB}$ is the $B$-orbit $C_{wB}:=Be_{w}$ in $G/B$.
Then the corresponding \red{\em Schubert variety}, denoted $X_{wB}$, is the Zariski closure $\overline{Be_{w}}$ in $G/B$. 
For an accessible introduction to the basic geometry of Schubert varieties, we recommend Brion's lectures~\cite{Brion2005}.

The partial order $\leq$ on $W$ defined by $v\leq w$ if and only if $X_{vB}\subseteq X_{wB}$, where $v,w\in W$, is called the \red{\em Bruhat-Chevalley order}. 
The complement $X_{wB}\setminus C_{wB}$ is called the \red{\em boundary of $X_{wB}$}, and denoted by $\partial X_{wB}$. 
We have
\begin{align*}
\partial X_{wB}=\bigcup\limits_{v\in W; v\lneq w}X_{vB}.
\end{align*}
For $v\,\in\, W,$ let $C^{vB}$ denote the orbit $B^{-}e_{v}$.
We call it the \red{\it opposite Schubert cell} in $G/B$. 
Then the \red{\it opposite Schubert variety corresponding to $v$} is defined as $\overline{B^{-}e_{v}}$.
To denote $\overline{B^{-}e_{v}}$ we write $X^{vB}$. 
Note that the complement $X^{vB}\setminus C^{vB}$ is the \red{\em boundary} of $\partial X^{vB}$. 
It is given by 
\begin{align*}
\partial X^{vB}=\bigcup\limits_{u\in W; v\lneq u}X^{uB}.
\end{align*}
Finally, we are ready to introduce the Richardson varieties properly. 
Let $v$ and $w$ be two elements from $W$. 
The corresponding \red{\em Richardson variety} in $G/B$, denoted $X_w^v$, is given by the intersection
\begin{align*}
X_w^v := X_{wB}\cap X^{vB}.
\end{align*} 
Schubert varieties are clearly Richardson varieties in the form $X_w^1$, where $1$ is the identity element of $W$.
This similarity extends to the properties of these varieties, as general Richardson varieties are also normal and Cohen-Macaulay, as proven in~\cite[Lemma 1]{BrionLakshmibai2003}.

A Richardson variety has two natural boundaries that are defined by 
\begin{align*}
(\partial X_{w})^{v}:= \partial X_{wB}\cap X^{vB}\quad \text{ and }\quad (\partial X^{v})_{w} := X_{wB}\cap \partial X^{vB}.
\end{align*}
We note in passing that $X_{w}^{v}$ and its boundaries are closed reduced, $T$-stable subvarieties of $G/B$. 
In~\cite[Corollary 1.5 and Theorem 2.1]{Richardson1992}, Richardson showed that the intersection $C_{wB}\bigcap C^{vB}$ is a smooth irreducible variety such that 
\begin{align*}
X_{w}^{v}=\overline{C_{wB}\bigcap C^{vB}}.
\end{align*}
Additionally, it follows readily from~\cite[Theorem 3.2]{Richardson1992} that the intersection $C_{wB}\bigcap C^{vB}$ is an open subset of $X_w^v$.
Moreover, it follows from~\cite[Corollary 1.2]{Deodhar1985} that $C_{wB}\cap C^{vB}$ is a nonempty open subset of $X_w^v$ if and only if $v\leq w$. 
By~\cite[Lemma 1 (i)]{BrionLakshmibai2003}, we know that $\dim X_w^v$ is given by $\ell(w)-\ell(v)$. 
We note here also that $C_{wB}$ is equal to the $U$-orbit $Ue_w$, where $U$ is the unipotent radical of $B$.
This follows from the fact that $B$ is a semidirect product of $U$ and $T$, and $e_w$ is a $T$-fixed point in $G/B$. 
Likewise, we have $C^{vB}=U^-e_v$, where $U^-$ is the unipotent radical of $B^-$.

Let $z$ and $y$ be two elements of $W$ such that $v\leq z$ and $v\leq y$.
Clearly, the relation $z\leq y$ in $W$ implies that $X_z^v\subseteq X_y^v$. 
Conversely, if the inclusion $X_z^v\subseteq X_y^v$ holds, then we have $e_z\in X_y^v$.
Since every $T$-fixed point $X_y^v$ is of the form $e_u$ for some $u$ from the interval $[v,y]$, we see that $z\leq y$. 
Hence, we proved the following simple lemma.
\begin{Lemma}\label{L:BConRichardson}
Let $v\leq w$. 
The interval $[v,w]$ is isomorphic to the inclusion order on the set of Richardson subvarieties $\{X_z^v\}_{v\leq z \leq w}$.
\end{Lemma}

We denote by $R$ the system of roots for the pair $(G,T)$.
The \red{\em set of positive roots} determined by $B$ will be denoted by $R^+$. 
Then the set of positive roots determined by the opposite Borel subgroup $B^{-}$, denoted $R^-$, is given by $R^{-} \,:=\, -R^{+}$.
The elements of $R^-$ (resp. $R$) are called the \red{\em negative roots} (resp. \red{\em positive roots}).
The indecomposable elements of $R^+$ are called the \red{\em simple roots}. 
We denote by $S$ the set of simple roots. 
Throughout this text, we use indexed Greek letters $\alpha_1,\alpha_2,\dots$ to denote simple roots in some root system $R$. 
If the set of simple roots is given by 
\begin{align}\label{A:ifSisgiven}
S=\{\alpha_1,\,\ldots,\,\alpha_n\}, 
\end{align}
then $n$ will be called the \red{\em rank of $G$}. 
For $S$ as in (\ref{A:ifSisgiven}), the simple reflection in $W$ corresponding to the simple root $\alpha_i$, where $1\leq i \leq n$, will be denoted by $s_{i}$.

The minimum of the number of simple reflections required to write an element $w\in W$ as a product is called the \red{\em length} of $w$,
and denoted by $\ell(w)$.
The  \red{\em length function}, $w\mapsto \ell(w)$, $w\in W$, is identical to the dimension function $w\mapsto \dim X_{wB}$, $w\in W$.
Furthermore, the length function serves as a grading on the Bruhat-Chevalley poset on $W$. 
A product of simple reflections $s_{i_1}s_{i_2}\cdots s_{i_r}$ is called a \red{\em reduced expression of $w$} if $w= s_{i_1}s_{i_2}\cdots s_{i_r}$ and $\ell(w) = r$. 
An element $c\in W$ is called a \red{\em Coxeter type element} if $c$ has a reduced expression $c=s_{i_1}s_{i_2}\cdots s_{i_r}$,
where $s_{i_1},\dots, s_{i_r}$ are mutually distinct simple reflections. 
\begin{Example}
Let $G$ denote the special linear group $SL(n,k)$.
The Weyl group of $G$, which is isomorphic to the symmetric group on $\{1,\dots, n\}$ will be denoted by $\mbf{S}_n$. 
The simple reflections that generate $\mbf{S}_n$ are given by $s_i$ ($1\leq i \leq n-1$), where $s_i$ is the permutation that interchanges $i$ and $i+1$ and leaves
everything else fixed. The Hasse diagram of the Bruhat-Chevalley order on $\mbf{S}_3$ is given in Figure~\ref{F:S3}:
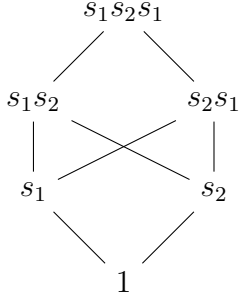
\begin{figure}
\begin{center}
	\begin{tikzpicture}[scale=.6]
		\node (e) at (0,0)  {$s_1s_2s_1$};
		\node (z) at (-2,-1.2) {};
		\node (y) at (2,-1.2) {};
		\node (g) at (-2,-2) {$s_1s_2$};
		\node (b) at (2,-2) {$s_2s_1$};
		\node (w) at (5.5,-4.8) {};
		\node (v) at (-2.5,-4.8) {};
		\node (u) at (2.5,-4.8) {};
		\node (d) at (-2,-4) {$s_1$};
		\node (p) at (2,-4) {$s_2$};
		\node (s) at (5.5,-7.6) {};
		\node (q) at (.5,-7.6) {};
		\node (a) at (0,-6) {$1$};
		\draw (e) -- (g) -- (d) -- (a) -- (p) -- (b) -- (e);
		\draw (b) -- (d);
		\draw (g) -- (p);
		\end{tikzpicture}
		\caption{The Bruhat-Chevalley order on $\mbf{S}_3$}
		\label{F:S3}
\end{center}
\end{figure}
\end{Example}

We will use the notation $X(T)$ (respectively, $Y(T)$) to denote \red{\em group of the characters} (respectively, the \red{\em group of co-characters}) of $T$. 
Both $X(T)$ and $Y(T)$ are free abelian groups.
They have equal rank. 
There is a bilinear perfect pairing between them, $\langle \cdot, \cdot \rangle: X(T)\times Y(T)\to \mathbb{Z}$.
Let ${X(T)}_{\mathbb{R}}:=X(T)\otimes _{\mathbb{Z}}\mathbb{R}$ and ${Y(T)}_{\mathbb{R}}:=Y(T)\otimes_{\mathbb{Z}}\mathbb{R}$. 
By a slight abusing of notation, we denote the induced bilinear pairing ${X(T)}_{\mathbb{R}}\times {Y(T)}_{\mathbb{R}}\to \mathbb{R}$ by $\langle \cdot, \cdot\rangle$, also. 
Let $(\cdot, \cdot)$ be a positive definite $W$-invariant symmetric bilinear form on ${X(T)}_{\mathbb{R}}$. 
For $\alpha \in R^{+}$, the reflection $s_{\alpha}$ can be defined by declaring 
\begin{align*}
s_{\alpha}(\lambda) :=\lambda-\frac{2(\lambda,\alpha)}{(\alpha,\alpha)}\alpha\qquad \text{for every $\lambda\in {X(T)}_{\mathbb{R}}$}.
\end{align*} 
For a root $\alpha\in R$, the associated \red{\em co-root} $\alpha^{\vee}\in {Y(T)}_{\mathbb{R}}$ is characterized by the property that for every $\lambda\in X(T)_\R$, 
the equality $\langle \lambda, \alpha^{\vee}\rangle=\frac{2(\lambda,\alpha)}{(\alpha,\alpha)}$ holds.
A \red{\em dominant character} is a character $\alpha \in X(T)$ such that the number $\langle \lambda, \alpha_{i}^{\vee}\rangle$ is a nonnegative integer for every $1\le i\le n$.
The set of dominant characters will be denoted by $X(T)_+$. 
For a simple root $\alpha_{i} \in S$, we denote the corresponding \red{\em fundamental weight} by $\varpi_{i}$.
This is the unique dominant weight with the property that $\langle \varpi_{i}, \alpha_{j}^{\vee}\rangle=\delta_{ij},$ where $\delta_{ij}$ denotes the Kronecker delta.

For each root $\alpha\in R$, there exists a 1-dimensional unipotent subgroup of $G$, which is normalized by $T$.
This specific 1-dimensional subgroup is denoted as $U_\alpha$, and is referred to as the \red{\em root subgroup associated with $\alpha$}. 
It is given by the image of a homomorphism $x_\alpha : (k,+)\to G$ satisfying the following identity 
for every $t\in k$ and $h\in T$: 
\begin{align*}
h x_\alpha(t) h^{-1} = x_\alpha(\alpha(h)t).
\end{align*}
It is worth noting that $G$ (resp. $B$) is generated by $T$ and all $U_{\pm \alpha}$ (resp. $T$ and all $U_\alpha$), where $\alpha \in S$.
It also worth noting that for every root $\beta\in R$, there is an element $w\in W$ and a simple root $\alpha\in S$ such that $\beta = w\alpha$
implying that the corresponding 1-dimensional unipotent subgroup $U_{\beta}$ is given by $U_\beta = w U_{\alpha} w^{-1}$.

\subsection{Subexpression of a reduced expression}

This subsection provides us with the main technical tool that is due to Deodhar~\cite{Deodhar1985}.
We will loosely follow the exposition of the article~\cite{MarshRietsch2004} of Marsh and Rietch.

Let $w\in W$.
An \red{\em expression} for $w$, denoted by $\bf{w}$, is a sequence $(w_{(0)},\dots,w_{(r)})\in W^{r+1}$ such that 
\begin{enumerate}
\item $w_{(0)}=1$,
\item $w_{(r)}=w$, 
\item $w_{(j)}$ is either $w_{(j-1)}$ or $w_{(j-1)}s_i$ for some $s_i\in S$. 
\end{enumerate}
Since for every $j\in \{1,\dots, r\}$, the product $w_{(j-1)}^{-1}w_{(j)}$ is an element of $S\cup \{1\}$, the information carried by the expression $\bf{w}$ is equivalent to its \red{\em sequence of factors}, that is, 
\begin{align*}
(w_{(1)} , w_{(1)}^{-1}w_{(2)}, \ldots, w_{(r-1)}^{-1}w_{(r)}).
\end{align*}

Now let us assume that the expression $\mbf{w} := (w_{(0)},\dots, w_{(r)})$ has length at least $2$. 
The following three subsets of indices will be convenient to use for some definitions: 
\begin{align*}
J_\mbf{w}^{+} & :=\{j\in \{1,\ldots, r\}: w_{(j-1)}<w_{(j)}\}, \\
J_\mbf{w}^{\circ} &:= \{j\in \{1,\ldots, r\}: w_{(j-1)}=w_{(j)}\}, \\
J_\mbf{w}^{-} &:= \{j\in \{1,\ldots, r\}: w_{(j)}<w_{(j-1)}\}.
\end{align*}
We say that an expression $\mbf{w}=(w_{(0)},w_{(1)},\ldots, w_{(r)})$ is \red{\em nondecreasing} (respectively, {\it \red{reduced}}) if 
$J_\mbf{w}^{-}=\emptyset$ (respectively, $J_\mbf{w}^- \cup J_\mbf{w}^\circ = \emptyset$) holds. 
Clearly, every nondecreasing expression $\mbf{w}$ gives rise to a reduced expression $\widehat{\mbf{w}}$ 
by discarding all $w_{(j)}$ from $\mbf{w}$ such that $j\in J_\mbf{w}^{\circ}$. 

\medskip

We are now in a position to define a delicate notion of a `positive subexpression.'
Let $(s_{i_1},s_{i_2},\ldots, s_{i_r})$ be the sequence of factors for a reduced expression $\mbf{w}$ for an element $w\in W$. 
Let $v\in W$ be an element such that $v\leq w$. 
A \red{\em subexpression for $v$ in $\mbf{w}$} is an expression $\mbf{v}:=(v_{(0)},v_{(1)},\ldots, v_{(r)})$ for $v$ such that 
\begin{align*}
v_{(j)}\in \{v_{(j-1)}, v_{(j-1)}s_{i_j}\}\quad\text{holds for every $j\in\{1,\ldots, r\}$.} 
\end{align*}
If, in addition, $v_{(j)}\leq v_{(j-1)}s_{i_{j}}$ holds for every $j\in \{1,\dots, r\}$, then we say that $\mbf{v}$ is a \red{\em distinguished subexpression for $v$ in $\mbf{w}$}.
In cases where the context implies the reference to the element $v\in W$, we will adopt a slight abuse of terminology by stating that `$\mathbf{v}$ is a subexpression of $\mathbf{w}$' or `$\mathbf{v}$ is a distinguished subexpression of $\mathbf{w}$.' Additionally, we will use the notation $\mathbf{v} \preceq \mathbf{w}$ to indicate that $\mathbf{v}$ is a distinguished subexpression of $\mathbf{w}$.

It is worth noting that within a distinguished subexpression $\mathbf{v} := (v_{(0)}, v_{(1)}, \ldots, v_{(r)})$, if there exists some $j \in \{1, \ldots, r\}$ such that the right multiplication by $s_{i_j}$ results in a decrease in the length of $v_{(j-1)}$, then we have $v_{(j)} = v_{(j-1)}s_{i_j}$.

Let $\mbf{w}$ be a reduced expression with sequence of factors $(s_{i_1},\dots, s_{i_r})$ as before. 
Let $\mbf{v}\preceq \mbf{w}$.
We will say that $\mbf{v}$ is a \red{\em positive subexpression} if the following (strict) relations hold: 
\begin{align}\label{A:positive}
v_{(j-1)} \lneq v_{(j-1)}s_{i_{j}}\qquad (j\in\{1,\ldots,r\}).
\end{align} 
We notice that the condition in (\ref{A:positive}) is equivalent to the condition that 
\begin{align*}
v_{(j-1)}\le  v_{(j)} \le v_{(j-1)}s_{i_j}\qquad (j\in\{1,\ldots,r\}).
\end{align*} 

From now on, we will use a `+' sign as a subscript for $\mathbf{v}$ to denote that it represents a positive subexpression.
We conclude our preparatory section by a remark.

\begin{Remark}\label{R:tworemarks}
\begin{enumerate}
\item It follows from definitions that $\mathbf{v}_+$ is characterized by the following two properties:
\begin{enumerate}
\item $\mathbf{v}_+\preceq \mathbf{w}$ (the distinguished subexpression property), 
\item $J_{\mbf{v}_+}^{-} =\emptyset$ (the nondecreasing property).
\end{enumerate}
\item For every $v\leq w$ in $W$ and a reduced expression $\mbf{w}$ of $w$, there is a unique positive subexpression $\mbf{v}_+$ of $v$ in $\mbf{w}$. 
A proof of this fact can be found in~\cite[Lemma 3.5]{MarshRietsch2004}.
\end{enumerate}
\end{Remark}

\begin{Example}\label{E:firstnontoric}
Let us consider for $G=SL(5,k)$ the elements $w:=s_3s_2s_1s_4s_3s_4s_2s_3$ and $v:=s_2s_1s_3s_4$.
Clearly, the product of the underscored entries in $s_3\underline{s}_2 \underline{s}_1s_4\underline{s}_3\underline{s}_4 s_2s_3$ gives $v$,
hence, $v\leq w$ in the Bruhat-Chevalley order.
A reduced expression for $w$ is given by 
\begin{align*}
\mbf{w} = (1,  s_3,   s_3s_2, s_3s_2s_1, s_3s_2s_1s_4,  s_3s_2s_1 s_4s_3,   s_3s_2s_1 s_4s_3 s_4,   s_3s_2s_1 s_4s_3 s_4s_2,  s_3s_2s_1 s_4s_3 s_4s_2s_3). 
\end{align*}
The sequence of factors for $\mbf{w}$ is given by $(s_3,s_2,s_1,s_4,s_3,s_4,s_2,s_3)$.
Here, we have two distinguished subexpressions for $v$ in $\mbf{w}$: 
\begin{align*}
\mbf{v}_1 &:= (1, 1, s_2 ,  s_2s_1,  s_2s_1 ,   s_2s_1 s_3, s_2s_1s_3s_4, s_2s_1s_3s_4,s_2s_1s_3s_4),\\
\mbf{v}_2 &:= (1, 1, s_2, s_2s_1,  s_2s_1s_4, s_2s_1s_4 s_3, s_2s_1s_4s_3s_4, s_2s_1s_4s_3s_4, s_2s_1s_3s_4).
\end{align*}
For the first distinguished subexpression we see that $J_{\mbf{v}_1}^\circ = \{ 1, 4,7,8\}$,  $J_{\mbf{v}_1}^+ = \{ 2,3,5,6\}$, $J_{\mbf{v}_1}^- = \emptyset$. 
The sequence of factors of $\mbf{v}_1$ is given by $(1,s_2,s_1,1,s_3,s_4,1,1)$.
Note that $\mbf{v}_1$ is nondecreasing. 
In other words, $\mbf{v}_1$ is the unique positive subexpression for $v$ in $\mbf{w}$.
For the second distinguished subexpression we have $J_{\mbf{v}_2}^\circ = \{ 1,7\}$, $J_{\mbf{v}_2}^+ = \{ 2,3,4,5,6\}$, $J_{\mbf{v}_2}^- = \{ 8\} $.
The sequence of factors of $\mbf{v}_2$ is given by $(1,s_2,s_1,s_4,s_3,s_4,1,s_3)$.
\end{Example}

\subsection{$R$-polynomials}
The \red{\em right descent set} of an element $w\in W$ is the set $D(w)$ of simple reflections defined by
\begin{align*}
D(w):= \{s_\alpha :\ \alpha \in S, \ \ell(ws_\alpha) < \ell(w) \}.
\end{align*} 
It is well-known from the work of Kazhdan and Lusztig~\cite{KazhdanLusztig1979} that there is a unique family of polynomials $\{R_{u,v}(q)\}_{u,v\in W} \subset \Z[q]$ satisfying the following conditions:
\begin{enumerate}
\item[(1)] $R_{u,v}(q) = 0$ if $u\nleq v$,
\item[(2)] $R_{u,v}(q) =1$ if $u=v$,
\item[(3)] If $s\in D(v)$, then 
$R_{u,v}(q) = 
\begin{cases}
R_{us,vs}(q) & \text{ if $s\in D(u)$},\\
qR_{us,vs}(q)+(q-1)R_{u,vs}(q) & \text{ if $s\notin D(u)$}.
\end{cases}$
\end{enumerate}
The polynomial $R_{u,v}(q)$ is referred to as the \red{\em $R$-polynomial} of the interval $[u,v]$.
Let $\F_q$ denote the finite field with $q$ elements, where $q$ is the power of a prime number of the form $q=p^m$, where $m\in \Z_+$.
For a variety $X$ defined over $\F_p$, we denote by $|X|_q$ the number $\F_q$-rational points of $X$. 
According to~\cite[Lemmas A3 and A4]{KazhdanLusztig1979}, $R$-polynomials have the following interpretation: 
\begin{align}\label{A:Ruv}
R_{u,v} (q):= | C^{uB}\cap C_{vB} |_q.
\end{align}
At the same time, for certain intervals, the corresponding $R$-polynomials take a special form. 
\begin{Lemma}\label{L:Brenti's} 
Let $[u,v]$ be an interval in $W$. 
Then the following statements are equivalent: 
\begin{enumerate}
\item[(i)] $[u,v]$ does not contain any interval isomorphic to the Bruhat-Chevalley order on $\mbf{S}_3$.
\item[(ii)] $R_{u,v}(q) = (q-1)^{\ell(v)-\ell(u)}$.
\end{enumerate}
\end{Lemma}
This lemma is proved by Brenti in~\cite[Theorem 6.3]{Brenti1994}. 
Hereafter, we will refer to this lemma as \red{\em Brenti's criterion}. 
Let us call an interval $[u,v]$ satisfying (i) an \red{\em $\mbf{S}_3$-free interval}.

\section{Torus actions on Richardson varieties}

\begin{Theorem}\label{T:arbitraryT}
Let $[v, w]$ be an interval in $W$. 
Let $\mc{T}$ be a $k$-torus. 
If the Richardson variety $X_{w}^{v}$ is a toric $\mc{T}$-variety, then all Richardson subvarieties $X_y^x\subseteq X_{w}^{v}$,
where $v\leq x\leq y\leq w$, are stable under the $\mc{T}$-action. 
\end{Theorem}

\begin{proof}
Since $X_w^v$ is a toric $\mc{T}$-variety, without loss of generality, we will proceed with the assumption that $\dim X_w^v = \dim \mc{T}$. 
Since $X_w^v$ is a projective variety, we know from Demazure's Structure theorem~\cite{Demazure1970} that $\Aut(X_{w}^{v})$ is a linear algebraic group.
We remark in passing that Demazure proved his result for the smooth toric varieties but in the article~\cite{LiendoArteche} Liendo and Arteche show that Demazure's proof remains true for every complete toric variety. 
Now, we know that $\mc{T}\subseteq \Aut(X_w^v)$ is a maximal torus for dimension reasons. 
Since $\mc{T}$ is connected, $\mc{T}$ is contained in the connected component of $\Aut(X_w^v)$ containing the identity element.
Let $\Aut^0(X_w^v)$ denote this identity component.  

Recall that $T$ denotes the maximal torus of $G$. 
Let $\varphi: T\to \Aut^0(X_{w}^{v})$ be the morphism of algebraic groups induced by the natural action of $T$ on $X_{w}^{v}$. 
Let $T'=\varphi(T).$ 
Without loss of generality, we may assume that $T'$ is contained in $\mc{T}$. 
We will analyze the action of $\mc{T}$ on the set of $T$-fixed points, that is, $(X_w^v)^T:=\{e_u :\ v\leq u \leq w\}$, in $X_w^v$.  
Since this is a finite set and $\mc{T}$ is connected, $\mc{T}$ acts trivially on $(X_w^v)^T$.
In other words, we have $(X_w^v)^T\subseteq (X_w^v)^{\mc{T}}$. 
Let $U$ (resp. $U^-$) denote the unipotent radical of $B$ (resp. of $B^-$). 
We will show that the intersection $Ue_{w}\cap U^{-}e_{v}$ is a $\mathcal{T}$-stable subset of $X_w^v$. 

Let $\xi$ (resp. $\tau$) be a point of $Ue_w\cap U^- e_v$ (resp. of $\mc{T}$). 
Let $\xi'$ denote the point obtained from the action of $\tau$ on $\xi$, that is, $\xi':= \tau \cdot \xi \in X_w^v$.
Since $\tau$ commutes with every element of $T$, $\tau$ gives a $T$-equivariant automorphism of $X_{w}^{v}$.
In particular, we have
\begin{equation*}
\overline{T\xi'}=\overline{\tau (T\xi)}= \tau (\overline{T\xi}).
\end{equation*}

By the Borel fixed point theorem, we know that every $T$-orbit closure in $G/B$ has a fixed point. 
In fact, every $B$-orbit in the flag variety contains a unique $T$-fixed point.
The unique $T$-fixed point in $B e_w$ is $e_w$. 
Since $B$ is given by the (semidirect) product $B=U\cdot T$, we see that $Ue_w = Be_w$.
In other words, we can characterize $Ue_w$ as the set of $x\in X_{wB}$ such that $\overline{T\cdot x}$ contains $e_w$ as its unique $T$-fixed point. 
Now recall our assumption that $\xi\in Ue_{w}$. 
Hence, $\overline{T\xi}$ contains $e_w$ as its unique fixed point. 
Therefore, we have 
\begin{equation*}
\overline{T\xi'}^{T}=\tau (\overline{T\xi}^{T})=\tau(e_{w})=e_{w}.
\end{equation*}
In words, the unique $T$-fixed point of the $T$-orbit closure of $\xi'$ is $e_w$.
Consequently, we see that $\xi'$ is an element of the orbit $Ue_w$. 
By using a similar reasoning, we see that $\xi'\in U^- e_v$. 
Hence, we have $\xi'\in Ue_w \cap U^{-}e_{v}$. 
It follows that $Ue_w \cap U^{-}e_{v}$ is $\mc{T}$-stable. 
Notice that this conclusion applies to every intersection of the form $Ue_y \cap U^-e_x$, where $x\leq y$ in $W$. 
But we know from Richardson's~\cite[Theorem 3.2]{Richardson1992} that $X_w^v = \bigsqcup_{v\leq x \leq y \leq w} U e_y \cap U^-e_x$. 
Hence, every Richardson subvariety $X_y^x:=\overline{Ue_y \cap U^-e_x}$ of $X_w^v$ is stable under the $\mc{T}$-action.
This finishes the proof of our theorem.	
\end{proof}	

\begin{Remark}
We expect that if a torus acts on a Richardson variety $X_w^v$ in $G/B$, then that torus must be isomorphic to a subtorus of the maximal torus $T$.
\end{Remark}

Another simple but useful corollary of our previous result is the following assertion.
\begin{Corollary}\label{C:S3isnottoric}
Let $X_w^v$ be a Richardson variety such that the interval $[v,w]$ is isomorphic to the Bruhat-Chevalley order on $\mbf{S}_3$.
Then $X_w^v$ is not a toric variety. 
\end{Corollary}

\begin{proof}
We assume towards a contradiction that the Richardson variety $X_w^v$ is a toric variety for some torus $\mathcal{T}$.
Let $\varphi : \mbf{S}_3\to [v,w]$ denote the poset isomorphism. 
Let $\{v=w_1,w_2,w_3,w_4,w_5,w_6=w\}$ be the set of elements of $[v,w]$. 
We list these elements so that we have $\varphi(1)=v,\varphi(s_1) = (w_2),\dots, \varphi(s_2s_1)=w_5$, and $\varphi(s_1s_2s_1)=w_6$. 
Then the Hasse diagram of the inclusion order on the Richardson subvarieties $\{X_z^v\}_{v\leq z\leq w}$, which is isomorphic to $\mbf{S}_3$, appears as in Figure~\ref{F:S3_second}. 
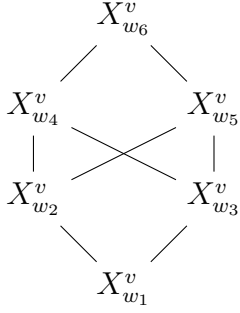
\begin{figure}[htp]
\begin{center}
	\begin{tikzpicture}[scale=.6]
		\node (e) at (0,0)  {$X_{w_6}^v$};
		\node (z) at (-2,-1.2) {};
		\node (y) at (2,-1.2) {};
		\node (g) at (-2,-2) {$X_{w_4}^v$};
		\node (b) at (2,-2) {$X_{w_5}^v$};
		\node (w) at (5.5,-4.8) {};
		\node (v) at (-2.5,-4.8) {};
		\node (u) at (2.5,-4.8) {};
		\node (d) at (-2,-4) {$X_{w_2}^v$};
		\node (p) at (2,-4) {$X_{w_3}^{v}$};
		\node (s) at (5.5,-7.6) {};
		\node (q) at (.5,-7.6) {};
		\node (a) at (0,-6) {$X_{w_1}^v$};
		\draw (e) -- (g) -- (d) -- (a) -- (p) -- (b) -- (e);
		\draw (b) -- (d);
		\draw (g) -- (p);
		\end{tikzpicture}
		\caption{The inclusion order on the Richardson varieties $\{X_z^v\}_{z\in [v,w]}$.}
		\label{F:S3_second}
\end{center}
\end{figure}
In fact, each of them is a $\mc{T}$-orbit closure. 
By Theorem~\ref{T:arbitraryT}, each of the Richardson varieties in $\{X_z^v\}_{z\in [v,w]}$ is a toric variety. 
In fact, each of them is a $\mc{T}$-orbit closure. 
In any toric variety, the torus orbit closures form a lattice with respect to inclusion. 
In particular, in our case, the $\mc{T}$-orbit closures in $X_w^v$ is a lattice. 
Hence, the intersection of any two $\mc{T}$-orbit closures have a unique maximal dimensional $\mc{T}$-orbit contained in them. 
Let us look at the intersection of the $\mc{T}$-orbit closures $X_{w_4}^{v}$ and $X_{w_5}^{v}$.
It is easily seen from Figure~\ref{F:S3_second} that 
\begin{align*}
X_{w_4}^v \cap X_{w_5}^v = X_{w_2}^v \cup X_{w_3}^v.
\end{align*}
But $X_{w_i}^v$ is isomorphic to $\PP^1$ for $i\in \{2,3\}$, and $X_{w_2}^v \cap X_{w_3}^v = X_{w_1}^v = e_v$.
Thus, the intersection of two $\mc{T}$-orbit closures, $X_{w_4}^v \cap X_{w_5}^v$, is not a $\mc{T}$-orbit closure. 
This contradiction proves that $X_w^v$ is not a toric variety. 
\end{proof}

Let us paraphrase our previous corollary.

\begin{Corollary}\label{C:iftoricthen}
If a Richardson variety $X_w^v$ in $G/B$ is a toric variety, then the interval $[v,w]$ is $\mbf{S}_3$-free.
\end{Corollary}

\begin{Example}
Let $[v,w]$ be an interval in a Weyl group $W$. 
We assume that $[v,w] \cong [1,s_1s_2s_1]$.
It follows from Corollary~\ref{C:iftoricthen} that $X_w^v$ is not a toric variety. 
To quantify the number of torus orbits in $X_w^v$, let us proceed with the assumption that $X_w^v$ is defined over a finite field $\F_q$.

In Table~\ref{tab:mytable}, we give the $R$-polynomials determined by the nontrivial subintervals of $[1,s_1s_2s_1]$. 
On the left column, we have the corresponding intersection in the Richardson variety. 
On the right column, we have the corresponding $R$-polynomials. 
\begin{table}[htp]
    \centering
    \begin{tabular}{|c|c|}
        \hline
        $C_{s_1s_2s_1B}\cap C^{1B}$ & $(q-1)^3 + q(q-1)$ \\
        \hline
        $C_{s_1s_2s_1B}\cap C^{s_1B}$ & $(q-1)^2$ \\
        \hline
       $C_{s_1s_2s_1B}\cap C^{s_2B}$ & $(q-1)^2$ \\
        \hline
       $C_{s_1s_2s_1B}\cap C^{s_1s_2B}$ & $(q-1)$ \\
        \hline
       $C_{s_1s_2s_1B}\cap C^{s_2s_1B}$ & $(q-1)$ \\
        \hline
       $C_{s_1s_2B}\cap C^{1B}$ & $(q-1)^2$ \\
        \hline
       $C_{s_1s_2B}\cap C^{s_1B}$ & $(q-1)$ \\
        \hline
       $C_{s_1s_2B}\cap C^{s_2B}$ & $(q-1)$ \\
        \hline
       $C_{s_2s_1B}\cap C^{1B}$ & $(q-1)^2$  \\
       \hline
       $C_{s_2s_1B}\cap C^{s_1B}$ & $(q-1)$ \\
        \hline
       $C_{s_2s_1B}\cap C^{s_21B}$ & $(q-1)$ \\
        \hline
       $C_{s_2B}\cap C^{1B}$ & $(q-1)$ \\
       \hline
       $C_{s_1B}\cap C^{1B}$ & $(q-1)$ \\ 
         \hline
    \end{tabular}
    \caption{The number of $\F_q$-rational points of the positive dimensional intersections. }
    \label{tab:mytable}
\end{table}
In addition to these positive dimensional $T$-varieties arising as intersections, we have 6 $T$-fixed points. 

We learn from Table~\ref{tab:mytable} two things: 
1) The total number of $\F_q$-rational points of the Richardson variety $X_w^v$ is given by 
\begin{align*}
|X_w^v|_q =  (q-1)^3 + 5(q-1)^2 + 9(q-1) + 6 = q^3+2q^2+2q +1.
\end{align*}
2) The maximal dimensional $T$-orbits in $C_{wB}\cap C^{vB}$ form a 1-parameter family over $\mathbb{A}^1$. 
Indeed, since $\dim T=2$, each maximal dimensional $T$-orbit in $C_{wB}\cap C^{vB}$ contributes $(q-1)^2$ rational points to $|C_{wB}\cap C^{vB}|_q$.
Since $|C_{wB}\cap C^{vB}|_q$ is given by $(q-1)^3+q(q-1)=q(q-1)^2  + (q-1)$,  we see that the number of $\F_q$-rational points on this family of maximal dimensional $T$-orbits is given by $q(q-1)^2$.
\end{Example}

\section{Characterizations of toric Richardson varieties}\label{S:Characterizations}
 
In this section, we give various characterization of the toric Richardson varieties for the natural action of the maximal torus $T$. 
We begin with an observation that relates the dimension of a Richardson variety to an appropriate reduced expression. 

\medskip

Let $X_w^v$ be a Richardson variety such that $\dim(X^{w}_{v})=l$.
Let $\mbf{w}$ be a reduced expression for $w\in W$ with sequence of factor $(s_{i_1},s_{i_2},\ldots, s_{i_r})$. 
Let $\mbf{v}_{+}=(v_{(0)},\ldots, v_{(r)})$ denote the unique positive subexpression for $v$ in $\mbf{w}$. 
Discarding from $\mbf{v}_+$ the entries ${v}_{(m)}$, where $m\in J_{\mbf{v}_{+}}^{\circ}$, gives us a new, reduced expression $\widehat{\mbf{v}}_{+}$ for $v$. 
This means that we have $\ell(v)=\ell(w)- |J_{\mbf{v}_{+}}^{\circ}|$.
In other words, we have the following relationship between the dimension of our Richardson variety and the set $J_{\mbf{v}_{+}}^{\circ}$:
\begin{align*}
\dim (X_{w}^{v})= |J_{\mbf{v}_{+}}^{\circ}|.
\end{align*}

To facilitate the proof of our next main theorem, we need additional notation and remarks. 
Let $w_1,v_1$, and $v_1'$ be elements from $W$ that satisfy $w_1 = v_1v_1'$ and $\ell(w_1) = \ell(v_1)+\ell(v'_1)$.
In particular, we have $v_1\leq w_1$. 
Since we have $w_1= v_1v'_1$ the existence of the following morphism is guaranteed: 
\begin{align*}
\pi^{w_1}_{v_1} : B e_{w_1} &\longrightarrow B e_{v_1} \\
b e_{w_1} &\longmapsto be_{v_1}.
\end{align*}
We will consider such morphisms in the context of distinguished subexpressions. 
To this end, let $\mbf{w}=(w_{(0)},\dots, w_{(r)})$ be a reduced expression for $w$ with the associated sequence of factors, $(s_{i_1},s_{i_2},\ldots, s_{i_r})$. 
Then for each $j\in \{1,\dots, r\}$, we have the corresponding morphism $\pi^w_{w(j)} : Be_w \to B e_{w(j)}$.
Now let $v\in W$ be an element such that $v\leq w$. 
For brevity, we will write $\mc{R}_{v,w}$ in place of $U e_{w} \cap U^-e_{v}$. 
Since $T$ fixes $e_{w}$ and $e_{v}$, and since $B=UT$ and $B^-=U^-T$, we have $\mc{R}_{v,w} = Be_{w}\cap B^- e_{v}$. 
For $x\in \mc{R}_{v,w}$ and $j\in \{1,\dots, r\}$, we denote by $x_j$ the image of $x$ under the morphism $\pi^w_{w(j)}$.  
By using the decomposition $G=\bigsqcup_{w\in W} B^-n_wB$, we will define a new element $v_{(j)}\in W$ such that $x_j \in B^-e_{v(j)}$. 
By the uniqueness of the decomposition $G=\bigsqcup_{w\in W} B^-n_wB$, we know that every element $x_j$ is contained in a unique coset of the form $B^- n_u B/B$ for some $u\in W$. 
We denote this $u$ by $v_{(j)}$, and let $\mbf{v}$ denote the sequence $\mbf{v}:=(v_{(0)},\dots, v_{(r)})$. 
Following the notation of~\cite[\S 4.1]{MarshRietsch2004}, we define 
\begin{align*}
\mc{R}_{\mbf{v},\mbf{w}} := \{x\in\mc{R}_{v,w} :\ \pi^w_{w(k)}(x) \in \mc{R}_{v_{(k)},w_{(k)}}\},
\end{align*}
and refer to it as a \red{\em Deodhar component}.
In~\cite[Theorem 1.1 and Corollary 1.2]{Deodhar1985}, Deodhar shows that 
\begin{enumerate}
\item $\mbf{v}\preceq \mbf{w}$ if and only if $\mc{R}_{\mbf{v},\mbf{w}}$ is nonempty, 
\item if $\mbf{v}\preceq \mbf{w}$ holds, then $\mc{R}_{\mbf{v},\mbf{w}} \cong (k^*)^{ | J_{\mbf{v}}^\circ|} \times k^{ | J_{\mbf{v}}^-|}$.
\end{enumerate}
It follows as a corollary of this theorem that $\mc{R}_{v,w} = \bigsqcup_{\mbf{v}\preceq \mbf{w}} \mc{R}_{\mbf{v},\mbf{w}}$,
where the union is taken over all distinguished subexpressions of $v$ in $\mbf{w}$. 
Clearly, this decomposition provides us with a further refinement of the Bruhat-Chevalley decomposition of $G$. 

We want to mention several other important properties of Deodhar components.  
\begin{Lemma}\label{L:stratum}
Let $\mbf{w}$ be a reduced expression for $w$. 
Let $v\leq w$. 
Then we have: 
\begin{enumerate}
\item[(i)] each Deodhar component in $\mc{R}_{v,w} = \bigsqcup_{\mbf{v}\preceq \mbf{w}} \mc{R}_{\mbf{v},\mbf{w}}$ is $T$-stable;
\item[(ii)] the Deodhar component $\mc{R}_{\mbf{v}_+,\mbf{w}}$ is dense in $\mc{R}_{v,w}$.
\end{enumerate} 
\end{Lemma}
The proof of (i) follows readily from the fact that the morphisms $\pi_{w(k)}^w$ are $T$-equivariant.  
The proof of (ii) follows from the fact that $|J_{\mbf{v}_+}^\circ | = \ell(w)-\ell(v)$. 
\medskip

We maintain our notation from above.
In addition, let $S=\{\alpha_1,\dots, \alpha_n\}$ denote the set of simple roots for our triplet $(G,B,T)$. 
Let $(s_{i_1},\dots, s_{i_r})$ be the sequence of factors of $\mbf{w}$. 
As before, we set $J_{\mbf{v}_+}^\circ:=\{j_1,\dots,j_l\}$, where $1\leq j_1<\cdots < j_l \leq r$.
Then, for $1\le m\le l$, we have 
\begin{align}\label{A:specialv}
v_{(j_m)} :=
s_{i_{1}}\cdots s_{i_{j_{1}-1}}\widehat{s}{_{i_{j_1}}}s_{i_{j_{1}+1}}\cdots s_{i_{j_{2}-1}}\widehat{s}{_{i_{j_2}}}\cdots \widehat{s}{_{i_{j_{m-1}}}}s_{i_{j_{m-1}+1}}\cdots s_{i_{j_{m}-1}}\widehat{s}{_{i_{j_{m}}}}. 
\end{align}
We are now ready to introduce a novel gadget.
We associate to the pair $(\mbf{v}_+,\mbf{w})$ a set of special roots by setting 
\begin{align}\label{A:specialroots}
\beta_{i_{j_m}}:={v_{(j_m)}}(\alpha_{i_{j_{m}}})\qquad  (1\le m\le l).
\end{align}

\begin{Remark}\label{distinguished}
Since $\mbf{v}_+$ is a positive subexpression for $v$ in $\mbf{w}$, by (\ref{A:positive}), the roots $\beta_{i_{j_m}}$ ($m\in \{1,\dots, l\}$) are positive roots. 
\end{Remark}

\begin{Example}\label{E:toricforSL5}
Let us consider for $G=SL(5,k)$ the elements $w:=s_3s_2s_1s_4s_3s_4s_2s_3$ and $v:=s_2s_1s_3s_2$.
Then we have $v\leq w$. 
A reduced expression for $w$ is given by 
\begin{align*}
\mbf{w} = (1,  s_3,   s_3s_2, s_3s_2s_1, s_3s_2s_1s_4,  s_3s_2s_1 s_4s_3,   s_3s_2s_1 s_4s_3 s_4,   s_3s_2s_1 s_4s_3 s_4s_2,  s_3s_2s_1 s_4s_3 s_4s_2s_3). 
\end{align*}
The sequence of factors for $\mbf{w}$ is given by $(s_3,s_2,s_1,s_4,s_3,s_4,s_2,s_3)$.
It is easily checked that there is a unique distinguished subexpressions for $v$ in $\mbf{w}$, that is, 
\begin{align*}
\mbf{v}_+ &:= (1, 1, s_2 ,  s_2s_1,  s_2s_1 ,   s_2s_1 s_3,       s_2s_1s_3,       s_2s_1s_3s_2,      s_2s_1s_3s_2).
\end{align*}
The sequence of factors for $\mbf{v}_+$ is given by $(1,s_2,s_1,1,s_3,1,s_2,1)$.
Then we have $J_{\mbf{v}_+}^\circ = \{ 1, 4, 6,8\}$.
Let $\{j_1,j_2,j_3,j_4\}$ denote the $J_{\mbf{v}_+}^\circ$.
Hence, we have 
\begin{align*}
v_{(j_1)} &= v_{(1)} = 1, \\
v_{(j_2)} &= v_{(4)} = s_2s_1, \\
v_{(j_3)} &= v_{(7)} = s_2s_1s_3, \\
v_{(j_4)} &= v_{(8)} = s_2s_1s_3s_2.
\end{align*}
It follows that 
\begin{align*}
\beta_{j_1} &= v_{(1)} (\alpha_3) = 1 (\alpha_3)=\alpha_3, \\
\beta_{j_2} &= v_{(4)} (\alpha_4) = s_2s_1(\alpha_4) = \alpha_4, \\
\beta_{j_3} &= v_{(6)} (\alpha_4) = s_2s_1s_3 (\alpha_4)= \alpha_2+ \alpha_3+\alpha_4, \\
\beta_{j_4} &= v_{(8)} (\alpha_3) = s_2s_1s_3s_2 (\alpha_3) =  \alpha_1.
\end{align*}
\end{Example}

We have one more example where we will show that the roots $\beta_{i_{j_1}},\dots, \beta_{i_{j_l}}$ may have dependencies. 
\begin{Example}\label{E:nontoricforSL5}
Let us consider for $G=SL(5,k)$ the elements $w:=s_3s_2s_1s_4s_3s_4s_2s_3$ and $v:=s_2s_1s_3s_4$.
Then $v$ is as in Example~\ref{E:firstnontoric}. 
We continue our calculation from there. 
Let $\{j_1,j_2,j_3,j_4\}$ denote the set $J_{\mbf{v}_+}^\circ = \{ 1, 4, 7,8\}$.
Hence, we have 
\begin{align*}
v_{(j_1)} &= v_{(1)} = 1, \\
v_{(j_2)} &= v_{(4)} = s_2s_1, \\
v_{(j_3)} &= v_{(7)} = s_2s_1s_3s_4, \\
v_{(j_4)} &= v_{(8)} = s_2s_1s_3s_4.
\end{align*}
It follows that 
\begin{align*}
\beta_{j_1} &= v_{(1)} (\alpha_3) = 1 (\alpha_3)=\alpha_3, \\
\beta_{j_2} &= v_{(4)} (\alpha_4) = s_2s_1(\alpha_4) = \alpha_4, \\
\beta_{j_3} &= v_{(7)} (\alpha_4) = s_2s_1s_3 s_4(\alpha_2)= \alpha_1+\alpha_2+\alpha_3, \\
\beta_{j_4} &= v_{(8)} (\alpha_3) = s_2s_1s_3s_4 (\alpha_3) =   \alpha_4.
\end{align*}
\end{Example}

The following assertion is one of the main results of our paper. 
\begin{Theorem}\label{th3.4}
Let $X_{w}^v$ be an $l$-dimensional Richardson variety. 
If $\{\beta_{i_{j_m}}: 1\le m\le l\}$ is the set of roots as defined in (\ref{A:specialroots}), 
then the Richardson variety $X_{w}^{v}$ is a toric $T$-variety if and only if the set $\{\beta_{i_{j_m}}: 1\le m\le l\}$ is linearly independent. 
\end{Theorem}

\begin{proof}
We know that $\mc{R}_{v,w}$ is an open subset of $X_w^v$.
It follows from Lemma~\ref{L:stratum} that $\mc{R}_{v,w}$ is a toric variety if and only if the dense stratum $\mathcal{R}_{\mbf{v}_{+}, \mbf{w}}$ is a toric variety.
Therefore, to prove our theorem, it suffices to show that the assertion of our theorem holds if and only if $\mathcal{R}_{\mbf{v}_{+}, \mbf{w}}$ is a toric variety.
\medskip

Following the notation of~\cite[\S 5.1]{MarshRietsch2004}, we define the following subset of $G$ for $\mbf{v}\preceq \mbf{w}$:
\begin{equation*}
G_{\mbf{v},\mbf{w}}:=
\left\{
g=g_1g_2\cdots g_r \Bigg\vert\
\begin{matrix}
g_j = x_{i_j} (t_j) s_{i_j} &\text{if $j\in J^-_{\mbf{v}}$} \\ 
g_j = y_{i_j} (a_j) &\text{if $j\in J^\circ_{\mbf{v}}$} \\
g_j = s_{i_j} &\text{if $j\in J^+_{\mbf{v}}$} \\
\end{matrix}
\quad \text{for $a_j \in k^*$ and $t_j\in k$}
\right\}.
\end{equation*}
Here, $x_{i_j}$ (resp. $y_{i_j}$) stands for the homomorphism $x_{\alpha_{i_j}} : k\to G$ (resp. $y_{\alpha_{i_j}} : k\to G$) that defines the 1-dimensional unipotent subgroup $U_{\alpha_{i_j}}$ (resp. $U_{-\alpha_{i_j}}$). We note in passing that $y_{\alpha_{i_j}} = x_{-\alpha_{i_j}}$.

To give an alternative proof of a result of Deodhar, in~\cite[Proposition 5.2(2)]{MarshRietsch2004}, Marsh and Rietsch notice that there is an isomorphism 
$f_{\mbf{v},\mbf{w}}: (k^*)^{ | J_{\mbf{v}}^\circ|} \times k^{ | J_{\mbf{v}}^-|}\to G_{\mbf{v},\mbf{w}}$, which is given by the obvious map that is obtained from the parameters $t_j\in k$ and $a_j\in k^*$ indicated in the definition of $G_{\mbf{v},\mbf{w}}$.
Composing this isomorphism with the canonical projection $\pi: G\to G/B$, they see that $(k^*)^{ | J_{\mbf{v}}^\circ|} \times k^{ | J_{\mbf{v}}^-|}$ (hence, $G_{\mbf{v},\mbf{w}}$) is isomorphic to the Deodhar component $\mc{R}_{\mbf{v},\mbf{w}}$. 
We restrict their isomorphism onto $\mc{R}_{\mbf{v}_+,\mbf{w}}$.
In this case, we already know that $J_{\mbf{v}_+}^-=\emptyset$.
Let $\{j_1,\dots, j_l\}$ be the set $J_{\mbf{v}_+}^\circ$.
Then the isomorphisms of Marsh and Rietsch are given by 
\begin{align*}
f_{\mbf{v}_+,\mbf{w}} :\ (k^*)^{ | J_{\mathbf{v}_+}^\circ|} &\longrightarrow G_{\mathbf{v}_+,\mathbf{w}} \\
(a_1,\dots, a_l) &\longmapsto 
s_{i_{1}}\cdots s_{i_{j_{1}-1}} y_{i_{j_1}} (a_1) s_{i_{j_{1}+1}}\cdots s_{i_{j_{2}-1}}
y_{i_{j_2}} (a_2) s_{i_{j_{2}+1}}\cdots s_{i_{j_{l}-1}} y_{i_{j_l}} (a_l) s_{i_{j_{l}+1}}
\cdots s_{i_r}
\end{align*}
and 
\begin{align*}
\pi |_{G_{\mbf{v}_+,\mbf{w}}} :  G_{\mathbf{v}_+,\mathbf{w}} &\longrightarrow R_{\mathbf{v}_+,\mathbf{w}} \\
f_{\mbf{v}_+,\mbf{w}} (a_1,\dots, a_l) &\longmapsto 
f_{\mbf{v}_+,\mbf{w}} (a_1,\dots, a_l)B.
\end{align*}
Now, we rewrite the initial segment $s_{i_{1}}\cdots s_{i_{j_{1}-1}} y_{i_{j_1}} (a_1)$ of $f_{\mbf{v}_+,\mbf{w}} (a_1,\dots, a_l)$ as follows:  
\begin{align*}
s_{i_{1}}\cdots s_{i_{j_{1}-1}} y_{i_{j_1}} (a_1) &= (s_{i_{1}}\cdots s_{i_{j_{1}-1}}) y_{i_{j_1}} (a_1) (s_{i_{j_1-1}}\cdots s_{i_{1}}) (s_{i_{1}}\cdots s_{i_{j_{1}-1}})\\
&= y_{\beta_{i_{j_1}}}(a_1) s_{i_{1}}\cdots s_{i_{j_{1}-1}} \qquad (\text{by using~(\ref{A:specialroots})}).
\end{align*}
Hence, we have 
\begin{align*}
f_{\mbf{v}_+,\mbf{w}} (a_1,\dots, a_l)= 
y_{\beta_{i_{j_1}}} (a_1) s_{i_{1}}\cdots s_{i_{j_{1}-1}}s_{i_{j_{1}+1}}\cdots s_{i_{j_{2}-1}}y_{i_{j_2}} (a_2) 
s_{i_{j_{2}+1}}\cdots s_{i_{j_{l}-1}} y_{i_{j_l}} (a_l) s_{i_{j_{l}+1}}
\cdots s_{i_r}.
\end{align*}
By repeatedly applying the same procedure, we reach at the following expression:
\begin{align*}
f_{\mbf{v}_+,\mbf{w}} (a_1,\dots, a_l)=   \left( \prod_{m=1}^l y_{\beta_{i_{j_m}}} (a_m) \right)
s_{i_{1}}\cdots s_{i_{j_{1}-1}}\widehat{s}{_{i_{j_1}}}s_{i_{j_{1}+1}}\cdots  s_{i_{j_{l}-1}}\widehat{s}{_{i_{j_{l}}}}s_{i_{j_l +1}}\cdots s_{i_r}.
\end{align*}
But earlier in this section we pointed out the fact that by discarding from $\mbf{v}_+$ the entries $v_{(m)}$,
where $m\in \{i_{j_1},\dots, i_{j_l}\}$, we obtain a new reduced expression $\widehat{\mbf{v}}_+$ for $v$. 
The last entry of this reduced expression is given by  
\begin{align*}
v= s_{i_{1}}\cdots s_{i_{j_{1}-1}}\widehat{s}{_{i_{j_1}}}s_{i_{j_{1}+1}}\cdots  s_{i_{j_{l}-1}}\widehat{s}{_{i_{j_{l}}}}s_{i_{j_l +1}}\cdots s_{i_r}.
\end{align*}
In other words, we have $f_{\mbf{v}_+,\mbf{w}} (a_1,\dots, a_l)=  \left( \prod_{m=1}^l y_{ \beta_{i_{j_m}}} (a_m) \right) v$.
It follows from this computation that 
\begin{align*}
\mc{R}_{\mbf{v}_+,\mbf{w}} &= \left\{  \prod_{m=1}^l y_{\beta_{i_{j_m}}} (a_m) v B \in G/B \ \bigg\vert\ (a_1,\dots, a_l)\in (k^*)^{ | J_{\mathbf{v}_+}^\circ|} \right\} \\
&=U^*_{-\beta_{i_{j_1}}} U^*_{-\beta_{i_{j_2}}}\cdots U^*_{-\beta_{i_{j_l}}}  e_v,
\end{align*}
where $U^*_{-\beta_{i_{j_m}}}$ denotes the non-identity elements of the 1-dimensional unipotent group $U_{-\beta_{i_{j_m}}}$.
The $T$-action on $U^*_{-\beta_{i_{j_1}}} U^*_{-\beta_{i_{j_2}}}\cdots U^*_{-\beta_{i_{j_l}}}  e_v$ is given by 
\begin{align}\label{A:theaction}
h\cdot   \prod_{m=1}^l y_{\beta_{i_{j_m}}} (a_m) v B =  \prod_{m=1}^l y_{\beta_{i_{j_m}}} (\beta_{i_{j_m}} (h) a_m) v B
\end{align}
for $h\in T$ and $(a_1,\dots, a_l)\in (k^*)^{|J_{\mbf{v}_+}^\circ|}$.
But we already know that $U^*_{-\beta_{i_{j_1}}} U^*_{-\beta_{i_{j_2}}}\cdots U^*_{-\beta_{i_{j_l}}}  e_v$ is isomorphic to $(k^*)^l$ via $f_{\mbf{v}_+,\mbf{w}}$,
and the action in (\ref{A:theaction}) transfers to the action 
\begin{align*}
h\cdot (a_1,\dots, a_l) = (\beta_{i_{j_1}} (h) a_1,\ldots, \beta_{i_{j_l}} (h) a_l )
\end{align*}
for $h\in T$ and $(a_1,\dots, a_l)\in (k^*)^{|J_{\mbf{v}_+}^\circ|}$.
Thus, if $l \leq \dim T$ holds, then the natural $T$-action on $\mc{R}_{\mbf{v}_+,\mbf{w}}$ defines defines a $T$-action on $(k^*)^l$.
We now proceed to check the dimension of the stabilizer group of a point $\xi \in \mathcal{R}_{\mbf{v}_{+}, \mbf{w}}$ in $T$. 
Since we have 
${\rm Stab}_{T}(\xi)= \bigcap\limits_{m=1}^{l}\ker(-\beta_{i_{j_{m}}}) \subset T$, 
the dimension of  the stabilizer group is $n-l$ if and only if the roots $\beta_{i_{j_1}},\dots, \beta_{i_{j_l}}$ are linearly independent. 
Said differently, $\dim T\cdot \xi = l$ if and only if the set of roots $\{\beta_{i_{j_1}},\dots, \beta_{i_{j_l}}\}$ is linearly independent. 
But our variety $X_w^v$ is $l$-dimensional. 
Hence, we conclude that $\{\beta_{i_{j_1}},\dots, \beta_{i_{j_l}}\}$ is linearly independent if and only if $X_w^v$ is a toric variety.
This finishes the proof of our theorem.
\end{proof}

Next, we state a corollary of our previous result that captures the main result of~\cite{Karuppuchamy2013}. 
Recall that a Coxeter type element is a product of distinct simple reflections.

\begin{Corollary}(Karuppuchamy~\cite{Karuppuchamy2013})
Assume that $v=1$, and the rest fo the hypothesis of Theorem~\ref{th3.4}.
Then the following statements are equivalent: 
\begin{enumerate}
\item[(1)] the Schubert variety $X_{wB}$ is a toric $T$-variety,
\item[(2)] the set $\{\alpha_{i_{j_m}} :\ 1\leq m \leq l\}$ is linearly independent,
\item[(3)] $w$ is a Coxeter element.
\end{enumerate}
\end{Corollary}
\begin{proof}
Let $\mbf{w}$ be a reduced expression for $w$. 
Since $v=1$, we know that the Richardson variety $X_w^v$ is the Schubert variety $X_{wB}$. 
Moreover, we know that $J_{\mbf{v}_+}^\circ = \{1,\dots, l\}$. 
It follows that $\{\beta_{i_{j_m}} :\ 1\leq m \leq l\}$ is the set of simple roots corresponding to the simple reflections that appear in the sequence of factors of $\mbf{w}$,
\[
\{\beta_{i_{j_m}} :\ 1\leq m \leq l\} = \{\alpha_{i_{j_m}} :\ 1\leq m \leq l\}.
\]
Now, by Theorem~\ref{th3.4}, $X_{wB}$ is a toric variety if and only if the set of roots $\{\beta_{i_{j_m}} :\ 1\leq m \leq l\}$ is linearly independent.
This proves the equivalence of (1) and (2). 
The equivalence of (2) and (3) follows at once from the fact that $w$ is a product of simple reflections that correspond to the elements of $\{\beta_{i_{j_m}} :\ 1\leq m \leq l\}$.
This finishes the proof of our assertion.
\end{proof}

Let us denote by $\mc{X}(R^+)$ the set of all reflections determined by the system of positive roots, 
\begin{align*}
\mathcal{X}(R^+):=\{s_{\beta}: \beta \in R^{+}\}.
\end{align*}
Since $W$ is generated by $\mathcal{X}(R^+)$, every element $w\in W$ can be written as a product 
\begin{align}\label{A:lfactors}
w=s_{\beta_1}\cdots s_{\beta_l}\quad \text{ for some $\beta_1,\dots ,\beta_l\in R^+$}.
\end{align} 
The smallest number of reflections required to write $w$ as in (\ref{A:lfactors}) is called the \red{\em absolute length} of $w$.
We will use the notation $\ell_a(w)$ to denote the absolute length of $w$.
We know from~\cite[Lemma 2]{Carter1972} that for every $w\in W$ the absolute length $\ell_a(w)$ is always less than or equal to $l$.
Now, by combining this new notion with the ideas used in the proof of Theorem~\ref{th3.4}, we can establish the previously promised assertion.

\begin{Proposition}\label{P:promised}
Let $[v,w]$ be an interval in $W$ such that $\ell(w)-\ell(v) \leq \dim T$. 
If the $R$-polynomial of a Richardson variety $X_w^v$ is given by $(q-1)^{\ell(w)-\ell(v)}$, then $X_w^v$ is a toric variety with respect to the natural action of $T$.
\end{Proposition}

\begin{proof}
We will closely follow the notation of the proof of Theorem~\ref{th3.4}.

As we mentioned earlier, in~\cite[Proposition 5.2]{MarshRietsch2004} Marsh and Rietsch give parametrizations for the Deodhar components of $\mc{R}_{v,w}$. 
Then Deodhar's description~\cite[Theorem 1.3]{Deodhar1985} of the $R$-polynomial of $[v,w]$ takes the following form: 
\begin{align}\label{A:R}
R_{v,w}(q) = \sum_{\mbf{v}\preceq \mbf{w}} q^{|J_{\mbf{v}}^-|} (q-1)^{|J_{\mbf{v}}^\circ|},
\end{align}
where the summation runs over all distinguished subexpressions of $v$ in $\mbf{w}$. 
To simplify our notation, let $\nu$ denote $\ell(w)-\ell(v)$. 
Now, let us assume that $(q-1)^\nu = R_{v,w}(q)$.
By collecting together the terms with equal power of $(q-1)$ on the right hand side of (\ref{A:R}), we obtain the following identity:
\begin{align}\label{A:orderofzero}
(q-1)^\nu = \sum_{d=0}^{\nu} f_d(q) (q-1)^d,
\end{align}
where $f_d(q) \in \N[q]$, and $\text{degree}(f_d(q)) \leq \nu-d$.
Of course, if $f_\nu (q) \neq 0$, then all other $f_d(q)$'s on the right hand side of (\ref{A:orderofzero}) are 0, implying that $f_\nu(q)=1$. 
Let us assume momentarily that $f_{\nu}(q)=0$.
In this case, we must have $f_0(q)=0$ as well. 
This follows from the fact that $f_0(q) \in \N[q]$, and all terms $f_d(q) (q-1)^d$, $d\in \{1,\dots, \nu\}$, are divisible by $q-1$. 
Now we know that the right hand side of (\ref{A:orderofzero}) is of the form $f_1(q)(q-1)+f_2(q)(q-1)^2+\cdots + f_{\nu-1}(q)(q-1)^{\nu-1}$.
After dividing both sides by $q-1$, we get an identity of the form 
\begin{align*}
(q-1)^{\nu-1} = \sum_{d=1}^{\nu} f_d(q) (q-1)^{d-1}.
\end{align*}
Then, by applying a similar argument, we see that $f_1(q)$ must be zero. 
We repeat this procedure to see that $f_0(q)=\cdots = f_\nu(q) = 0$, which is absurd. 
Therefore, we conclude that $f_\nu(q) = 1$. 
An important consequence of this observation is that there exists a unique distinguished subexpression $\mbf{v}$ in $\mbf{w}$, implying that 
\begin{align*}
(k^*)^{|J_{\mbf{v}_+}^\circ|} \cong \mc{R}_{\mbf{v}_+,\mbf{w}} = \mc{R}_{v,w}. 
\end{align*}
Then it follows from the last part of the proof of Theorem~\ref{th3.4} that $\mc{R}_{v,w}$ is given by 
$\mc{R}_{v,w}=U^*_{-\beta_{i_{j_1}}} U^*_{-\beta_{i_{j_2}}}\cdots U^*_{-\beta_{i_{j_l}}}  e_v$,
where $U^*_{-\beta_{i_{j_m}}}$ denotes the set of non-identity elements of $U_{-\beta_{i_{j_m}}}$,
and $\beta_{i_{j_1}},\dots, \beta_{i_{j_l}}$ are the roots defined as in (\ref{A:specialroots}).
We set $\Gamma:=\{\beta_{i_{j_1}},\dots, \beta_{i_{j_l}}\}$.

We recall from Remark~\ref{distinguished} that $\Gamma$ consists of positive roots only. 
Let us denote by $\mc{X}(\Gamma)$ the corresponding set of reflections, 
\begin{align*}
\mc{X}(\Gamma) :=\{s_{\beta_{i_{j_{m}}}}: 1\le m\le l\}.
\end{align*} 
Let $W'$ denote the reflection subgroup generated by $\mc{X}(\Gamma)$ in $W$. 
We proceed with a basic observation regarding the indices of the elements of $\mc{X}(\Gamma)$. 
The fact that $\mbf{v}_+$ is the unique distinguished expression for $v$ in $\mbf{w}$ implies that $\beta_{j_{k_1}}\neq \beta_{j_{k_2}}$ for every $\{ k_1, k_2\}\subseteq \{1,\dots,l\}$ such that $k_1\neq k_2$. 
Indeed, towards a contradiction, let us assume otherwise that there exist $1\le k_1<k_2\le l$ such that $\beta_{i_{j_{k_1}}}=\beta_{i_{j_{k_2}}}$.
Then we define an expression $\mbf{v'}$ having the following properties: 
$J_{\mbf{v}'}^+=J_{\mbf{v}_+}^{+}\cup\{k_1\},$ $J_{\mbf{v}'}^\circ=J_{\mbf{v}_+}^{\circ}\setminus \{k_1,k_2\},$ and $J_{\mbf{v}'}^-=\{k_2\}$. 
It is easy to see such an expression $\mbf{v}'$ is a distinguished subexpression for $v$ in $\mbf{w}$ that is different from $\mbf{v}_+$, which is absurd.
Now, since  $\beta_{i_{j_{k_{1}}}}\neq \beta_{i_{j_{k_{2}}}}$ for $1\le k_1\neq k_2\le l$, by~\cite[Theorem 4.4]{Dyer1990}, we see that the pair $(W', R')$ is a Coxeter system. 
It follows that the absolute length of the product of all generators is equal to the number of generators,
\begin{align*}
\ell_a(s_{\beta_{i_{j_{l}}}}\cdots s_{\beta_{i_{j_{1}}}})=l.
\end{align*}
Then it follows from~\cite[Lemma 3]{Carter1972} that $\Gamma$ is a linearly independent set of positive roots. 
Hence, our Theorem \ref{th3.4} shows that $X_w^v$ is a toric variety with respect to the natural action of $T$.
This finishes the proof of our proposition.
\end{proof}

We have several remarks in order.

\begin{Remark}
In~\cite[Proposition 6.1]{Deodhar1985}, by using `Deodhar components,' Deodhar finds an EL-labeling\footnote{Note: in the earlier days of the EL-labeling concept, the terminology ``$L$-labeling'' was more commonly used. 
There is also a notion of a CL-labeling which is weaker than the notion of an EL-labeling. 
Although these two notions seem to have similar topological consequences, they are not the same as shown by Tiansi Li in his thesis,~\cite{TiansiLi}.} of the Bruhat-Chevalley order on $[v,w]$.
In our case, it can be seen without much difficulty that the indices of our positive roots $\beta_{i_{j_1}},\dots, \beta_{i_{j_l}}$ give the EL-labels
along the unique increasing chain that Deodhar considered in his~\cite[Proposition 6.1 (i)]{Deodhar1985}. 
All of this is related to the polytope $\mc{Q}_{v,w}$ whose 1-skeleton is the Bruhat graph of $[v,w]$.
The positive roots $\beta_{i_{j_1}},\dots, \beta_{i_{j_l}}$ give the labels of the directed path on $\mc{Q}_{v,w}$ that starts at $v$ and ends at $w$. 
\end{Remark}

Proposition~\ref{P:promised} shows that the interval considered in Example~\ref{E:toricforSL5} is $\mbf{S}_3$-free but the one in 
Example~\ref{E:nontoricforSL5} contains a subinterval that is isomorphic to $\mbf{S}_3$.

\begin{Remark}
In Proposition~\ref{P:promised}, our assumption that $\ell(w)-\ell(v) \leq \dim T$ cannot be omitted. 
Indeed, the $R$-polynomial of the interval $[s_2,s_3 s_2 s_1 s_2 s_3]$ in $\mbf{S}_4$ is given by $(q-1)^4$. 
However, the maximal torus of $SL(4,k)$ is 3-dimensional. 
In other words, the Richardson variety $X_{s_3 s_2 s_1 s_2 s_3}^{s_2}$ is not a toric variety for the natural action of $T$. 
\end{Remark}

Let $X_w^v$ be a Richardson variety. 
Let $\mbf{w}=(w_{(0)},w_{(1)},\ldots, w_{(r)})$ be a reduced expression for $w$.
Let $(s_{i_{1}},\ldots ,s_{i_{r}})$ be the sequence of factors of $\mbf{w}$. 
Let $\mbf{v}_+:=(v_{(0)}, v_{(1)},\ldots, v_{(r)})$ denote the unique positive subexpression of $v$ in $\mbf{w}$. 
As before, let $J_{\mbf{v}_+}^\circ$ be given by $\{j_1,\dots, j_l\}$. 
In our next example, we will show that even if there are coincidences among the simple reflections $s_{i_{j_1}},\dots, s_{i_{j_l}}$ 
the Richardson variety $X_{w}^{v}$ might still be a toric variety.

\begin{Example}
Let $G=SL(3,k)$.
Let $T$ denote the maximal diagonal torus in $G$. 
We will consider $X_{w}^v$, where $w=s_{1}s_{2}s_{1}$ and $v=1s_{2}1$. 
Then $\mbf{w}=(1,s_1,s_1s_2,s_1s_2s_1)$ is a reduced expression for $w$, 
and $\mbf{v}_{+}=(1,1,s_{2},s_{2})$ is a positive subexpression for $v$ in $\mbf{w}$. 
Then we see that $J_{{\mbf{v}_{+}}}^{\circ}=\{1,3\}$. 
Note that $s_{\alpha_{i_{1}}}=s_{\alpha_{i_{3}}}$, but, by Theorem~\ref{th3.4} we know that $X_{w}^{v}$ is a toric variety. 
Indeed, $\mathcal{R}_{\mbf{v}_{+},\mbf{w}}\simeq U^{*}_{-\alpha_{1}}s_{2}U^{*}_{-\alpha_{1}}B=U^{*}_{-\alpha_{1}}U^{*}_{-\alpha_{1}-\alpha_{2}}vB$,
where $*$'s on the root subgroups indicate that we consider only the non-identity elements from them. 
Since $\{\alpha_1,\alpha_1+\alpha_2\}$ is a linearly independent set of roots, we see that $X_{w}^{v}$ is a toric variety for the action of the maximal torus of $SL(3,k)$.
\end{Example}

\begin{Remark}
Let $w\in W$. Then the $R$-polynomial of the lower interval $[1,w]$ in $W$ is given by $(q-1)^{\ell(w)}$ if and only if $w$ is a Coxeter element.
\begin{proof}
Let $\mbf{w}$ be a reduced expression for $w$. 
Evidently, there is a unique distinguished subexpression for 1 in $\mbf{w}$. 
If $w$ is a Coxeter element, then by~\cite[Theorem 1.3]{Deodhar1985}, we have $R_{1,w}(q) = (q-1)^{\ell(w)}$. 
Conversely, if $R_{1,w}(q)$ is given by $(q-1)^{\ell(w)}$, then, by~\cite[Ch 5. Exercise 18]{BjornerBrenti}, we know that $[1,w]$ is a Boolean lattice. 
Now it follows from the Subword Property of the Bruhat-Chevalley order~\cite[Theorem 2.2.2]{BjornerBrenti} that there are $\ell(w)$ distinct simple reflections in a reduced expression of $w$. These simple reflections appear as the atoms of the lattice $[1,w]$. 
Hence, $w$ is a product of distinct simple reflections.
In other words, $w$ is a Coxeter element. 
\end{proof}
\end{Remark}

We are now ready to prove Theorem~\ref{intro:T1} from Introduction. 
We recall its statement for convenience.

\medskip
{\it 
%\begin{Theorem}\label{intro:T1}
Let $v$ and $w$ be two elements from $W$ such that $v\leq w$. 
If the length of the interval $[v,w]$ is bounded by the dimension of $T$, then the following assertions are equivalent:
\begin{enumerate}
\item[(1)] the interval $[v,w]$ is $\mbf{S}_3$-free;
\item[(2)] the interval $[v,w]$ is a lattice with respect to the Bruhat-Chevalley order;
\item[(3)] the Richardson variety $X_w^v$ is a toric variety. 
\end{enumerate}}
%\end{Theorem}

\begin{proof}[Proof of Theorem~\ref{intro:T1}]

(1)$\Rightarrow$(3)
By Brenti's criterion, we know that the $R$-polynomial of $X_w^v$ is given by $(q-1)^{\ell(w)-\ell(v)}$.
The rest of the proof follows from Proposition~\ref{P:promised}. 
\medskip

(3)$\Rightarrow$(1)
This implication follows from Corollary~\ref{C:S3isnottoric}.
\medskip

(1)$\Rightarrow$(2).
In the proof of Proposition~\ref{P:promised} we show that if $R_{v,w}(q)=(q-1)^{\ell(w)-\ell(v)}$, then $C_{wB}\cap C^{vB}$ is the open $T$-orbit in $X_w^v$.
Thanks to Theorem~\ref{T:arbitraryT}, we may apply this fact to every Richardson subvariety $X_y^x$ ($v\leq x\leq y\leq w$), 
proving that the intersection $C_{yB}\cap C^{xB}$ is a $T$-orbit. 
In particular, for two Richardson subvarieties $X_y^v$ and $X_z^v$, where $v\leq y,z\leq w$, the inclusion $X_y^v \subseteq X_z^v$ is an inclusion of $T$-orbit closures.
Therefore, by Lemma~\ref{L:BConRichardson}, $[v,w]$ is isomorphic to the inclusion order of torus orbit closures.
But it is well-known that the inclusion order on the set of torus orbit closures in a projective toric variety is anti-isomorphic to the face-lattice of a polytope. 
This finishes the proof of (1)$\Rightarrow$(2).
\medskip

Finally, we will prove the implication (2)$\Rightarrow$(1). 
Let us assume that the interval $[v,w]$ is a lattice but not $\mbf{S}_3$-free.
Notice that the Bruhat-Chevalley order on $\mbf{S}_3$ is not a lattice. 
It follows that no lattice can contain a subinterval that is isomorphic to the Bruhat-Chevalley order on $\mbf{S}_3$. 
This finishes the proof of our theorem.
\end{proof}

\begin{Example}\label{E:anexamplefromSL4}
Let $G=SL(4,k)$. We consider the elements $w:=s_{1}s_{2}s_{3}s_{4}$ and $v:=s_{2}s_{4}$. 
Since the interval $[v,w]$ is $\mbf{S}_3$-free, by Theorem~\ref{intro:T1}, $X_w^v$ is a toric Richardson variety. 
\end{Example}

We will apply Theorem~\ref{th3.4} to obtain the proof of the following practical characterization of certain toric Richardson varieties. 

\begin{Corollary}\label{Cor3.4}
Let $v, w\in W$ be such that $w=cv$ and $\ell(w)=\ell(c)+\ell(v)$, where $c$ is a Coxeter type element. 
Then $X_{w}^{v}$ is a toric variety for the natural action of $T$.
\end{Corollary}
\begin{proof}
If $w$ is of the form $cv$, where $c$ is a Coxeter type element and $\ell(w) = \ell(c) + \ell(v)$, then we have a reduced expression of $w$ such that 
\begin{align*}
w = \underbrace{ s_{i_1}s_{i_2}\cdots s_{i_l}}_{c} \underbrace{s_{i_{l+1}}s_{i_{l+2}}\cdots s_{i_r}}_{v}.
\end{align*}
Then, for every $m\in\{1,\dots,l\}$, the element $v_{(m)}\in W$ defined in (\ref{A:specialv}) is the identity element. 
Hence, our special roots defined in (\ref{A:specialroots}) are given by
\begin{align*}
\beta_{i_m}:=\alpha_{i_m}\qquad \text{for $1\le m\le l$}.
\end{align*}
Since these are all distinct simple roots, they form a linearly independent set. 
Hence, our corollary follows from Theorem~\ref{th3.4}.
\end{proof}

In the following example we will show that there are toric Richardson varieties $X_w^v$ where $w$ need not be of the form $u_{1}vu_{2}$,
where $u_i$ ($i\in \{1,2\}$) is a Coxeter type element.

\begin{Example}\label{E:anexamplefromSL4_II}
We continue our Example~\ref{E:anexamplefromSL4}.
Recall that $w=s_{1}s_{2}s_{3}s_{4}$ and $v=s_{2}s_{4}$. 
It is easily seen that $w$ is not of the form $u_1vu_2$, where neither $u_1$ nor $u_2$ is a Coxeter type element. 
\end{Example}

\section{Toric Richardson varieties in Grassmannians}\label{S:Grassmanns}

In this section, we will consider Richardson varieties in Grassmannians.  
We maintain our assumptions on our algebraic group, $G$. 
As before, we will work with the root system $R$, the set of simple roots $S:=\{\alpha_1,\dots,\alpha_n\}$, and the weight lattice $X(T)$ determined by $(G,B,T)$. 
Recall that a fundamental weight $\varpi_m$ ($1\leq m\leq n$) is an element of $X(T)$ such that $\langle \varpi_m, \alpha_i\rangle =\delta_{im}$ for every $\alpha_i\in S$. 
A fundamental weight $\varpi_{r}$ is said to be \red{\it minuscule} if $\varpi_{r}$ satisfies $\langle\varpi_{r}, \beta \rangle \le 1$ for all $\beta \in R^{+}$.
A maximal parabolic subgroup $P$ is called \red{\em minuscule} if the Weyl group translates of a highest weight vector span the space of global sections on $G/P$ of the ample generator of the Picard group of $G/P$. 
The following table, which can be extracted from~\cite[pg. 247, Exercise 15]{Bourbaki4-6} or~\cite[Remark 2.3]{LMSIV}, gives the complete list of minuscule weights for simple groups:
\vspace{.1cm}
\begin{center}
	\begin{tabular}{ |p{1.3cm}|p{3.3cm}|p{3.5cm}| }
		\hline
		\multicolumn{3}{|c|}{ List of Minuscule weights} \\
		\hline
		no.&Root system type& Minuscule weights \\
		\hline
		$1.$ &	$A_{n}$  & $\varpi_{1}, \varpi_{2}, ..., \varpi_{n}$      \\
		
		$2.$ &	$B_{n}$  & $\varpi_{n}$      \\
		 
		$3.$ &	$C_{n}$  & $\varpi_{1}$  \\
		
		$4.$ & $D_{n} $  & $\varpi_{1}, \varpi_{n-1}, \varpi_{n}$       \\
		
		$5.$ &	$E_{6}$  & $\varpi_{1},\varpi_{6}$       \\
		
		$6.$ &	$E_{7}$  &    $\varpi_{7}$    \\
		
		$7.$ & $E_{8}$  & none      \\
		
		$8.$ & $F_{4}$  & none      \\
			
		$9.$ & $G_{2}$  & none      \\
		\hline
	\end{tabular}
\end{center}

Let $J$ be a subset of $S$. 
The subgroup of $W$ that is generated by the simple reflections $s_\alpha$, $\alpha \in J$ will be denoted by $W_J$. 
The set of minimal length left coset representatives of $W_J$ in $W$ will be denoted by $W^J$, and will be called the \red{\em parabolic quotient of $W$ determined by $J$}. 
In~\cite{Proctor1984}, Proctor classified the parabolic quotients of Weyl groups which are lattices under the Bruhat-Chevalley order. 
In particular, he showed that aside from minuscule cases, there is one additional situation, when the root system is of type $G_2$. 
The following result, which is a specialization of~\cite[Theorem 7.1]{Stembridge1996} to Weyl groups, provides us with a detailed description of the structure of these posets. 
\begin{Lemma}\label{L:Stembridge}
Let $W$ be a Weyl group. 
Let $J$ be a proper subset of $S$. 
Then the following assertions are equivalent: 
\begin{enumerate}
\item $(W^J,\leq)$ is a lattice.
\item $(W^J,\leq)$ is a distributive lattice.
\item $(W^J,\leq_L)$ is a distributive lattice.
\item $(W^J,\leq) \cong (W^J,\leq_L)$. 
\end{enumerate}
\end{Lemma}
In this lemma, the notation $\leq_L$ stands for the \red{\em left weak order}, which is defined as the transitive closure of the following cover relations:
\begin{align*}
\text{$v$ is covered by $w$} \iff w= s v\quad \text{for some $s \in S$},
\end{align*}
where $v,w\in W$.

Let $J$ be a proper subset of $S$. 
Let $P$ denote the parabolic subgroup generated by $B$ and all the 1-dimensional root subgroups $U_{-\alpha}$, where $\alpha \in J$. 
For $v$ and $w$ from $W^J$, the corresponding \red{\em parabolic Richardson variety} is defined as the following intersection in $G/P$: 
\begin{align*}
X_{wP}^{vP} := \overline{BwP/P}\cap \overline{B^-vP/P}.
\end{align*}
It is easy to see that $X_{wP}^{vP}$ is nonempty if and only if $v\leq w$ in $(W^J,\leq)$. 
We are now ready to characterize the toric parabolic Richardson varieties in our special Grassmann varieties.

\begin{Theorem}\label{T:Coxetertypeinterval}
Let $J\subseteq S$ be a set of simple roots such that $(W^J,\leq)$ is a lattice. 
Let $P$ denote the parabolic subgroup determined by $J$. 
Then $X_{wP}^{vP}$ is a toric Richardson variety in $G/P$ with respect to the action of $T$ if and only if $w$ is of the form $cv$ for some Coxeter type element $c$ such that $\ell(w) = \ell(c)+\ell(v)$.
\end{Theorem}

\begin{proof}
Since $v\leq_L w$, there exists an element $u\in W$ such that $uv = w$ and $\ell(u)+\ell(v) = \ell(w)$. 

We will show that $u$ is a Coxeter type element if and only if $X_{wP}^{vP}$ is a toric variety. 
We fix a reduced expression $\mbf{u}=(u_{(0)},\ldots ,u_{(l)})$ for $u$.
Let $(s_{i_{1}},\ldots, s_{i_{l}})$ be the sequence of factors of $\mbf{u}$. 
%Then there exists $1\le k_1<k_{2}\le l$ such that $i_{k_1}=i_{k_2}.$ 
Likewise, let $\mbf{v}=(v_{(0)},v_{(1)},\ldots, v_{(t)})$ be a reduced expression for $v$, and let $(s_{i_{l+1}},\ldots, s_{i_r})$ be the sequence of factors of $\mbf{v}$,
where $r=l+t$. 
It is easy to check that 
\begin{align*}
\mbf{w}=(u_{(0)},u_{(1)},\ldots,u_{(l)},u_{(l)}v_{(1)},\ldots, u_{(l)}v_{(t)}) 
\end{align*} 
is a reduced expression for $w$, and that $(s_{i_{1}},\ldots, s_{i_{l}}, s_{i_{l+1}}\ldots s_{i_r})$ is the sequence of factors for $\mbf{w}$. 
Also, it is easy to check that the subexpression 
\begin{align*}
\mbf{v}_{+}=(\underbrace{1,1,\ldots, 1}_{l-times},{ v}_{(1)},\ldots , { v}_{(t)})
\end{align*} 
for $v$ in $\mbf{w}$ is the unique positive subexpression for $v$ in $\mbf{w}$. 
It follows that the special roots defined in (\ref{A:specialroots}) are the simple roots, 
\begin{align*}
\beta_{i_m}:=\alpha_{i_m}\qquad \text{for $1\le m\le l$}.
\end{align*}
Hence, we have the identification $\mc{R}_{\mbf{v}_+,\mbf{w}}  = U^*_{-\alpha_{i_1}} U^*_{-\alpha_{i_2}}\cdots U^*_{-\alpha_{i_l}} vB/B$. 
Since $v\in W^J$ and $\ell(w) =\ell(u)+\ell(v)$, we see that the projection of $\mc{R}_{\mbf{v}_+,\mbf{w}}$ into $G/P$ is an isomorphism.
But we know $\{\beta_{i_1},\dots, \beta_{i_l}\}$ is a linearly independent set of simple roots if and only if 
$U^*_{-\alpha_{i_1}} U^*_{-\alpha_{i_2}}\cdots U^*_{-\alpha_{i_l}} vB/B$ has an open $T$-orbit. 
Equivalently, $U^*_{-\alpha_{i_1}} U^*_{-\alpha_{i_2}}\cdots U^*_{-\alpha_{i_l}} vP/P$ has an open $T$-orbit if and only if 
$\{\beta_{i_1},\dots, \beta_{i_l}\}$ is a linearly independent set of simple roots.
In other words, $X_{wP}^{vP}$ is a toric variety with respect to the action of $T$ if and only if the simple reflections $s_{i_1},\dots, s_{i_l}$ are mutually distinct.
It follows that $X_{wP}^{vP}$ is a toric variety with respect to the action of $T$ if and only if $u= s_{i_1}\cdots s_{i_l}$ is a Coxeter type element. 
This finishes the proof of our theorem.
\end{proof}

Paraphrasing our previous result gives us a proof of our Theorem~\ref{intro:T2} from the introduction section. 

\begin{proof}[Proof of Theorem~\ref{intro:T2}.]
The list of Grassmannians covered by the hypothesis of our Theorem~\ref{intro:T2} is precisely the list of Grassmannians covered by the hypothesis of our Theorem~\ref{T:Coxetertypeinterval}.
Hence, the proof of our theorem follows from Theorem~\ref{T:Coxetertypeinterval}.
\end{proof}

\section{Rationally Smooth Toric Richardson Varieties}\label{S:Concluding}

Let $X$ be a complex algebraic variety.
Let $H^*(X,\Q)$ denote the total cohomology space of $X$ with rational coefficients. 
For $x\in X$, we will use the shorthand $H^*_x(X)$ for the local cohomology $H^*(X, X\setminus \{x\})$.
Let $\dim_x X$ denote the dimension of the local ring of $X$ at $x$. 
The variety $X$ is said to be \red{\em rationally smooth at $x$} if for all $y$ in a neighborhood of $x$ in the complex topology,
$H_y^{2\dim_x(X)}(X)$ is isomorphic to $\Q$, and $H_y^i (X) = 0$ for $i\neq 2\dim_x(X)$. 
The set of rationally smooth points of $X$ is open in complex topology, and contains all smooth points.  
In~\cite[Theorem 4.3]{KazhdanLusztig1980}, Kazhdan and Lusztig showed that the Poincar\'e polynomial of the local intersection homology group $IH^*_v(X_{wB},\Q)$ is given by the \red{\em Kazhdan-Lusztig polynomial}, denoted $P_{v,w}(q)$. 
These polynomials can be defined axiomatically as follows: 

There exists a unique family of polynomials $\{P_{u,v}(q)\}_{u,v\in W} \subseteq \Z[q]$ such that for every $u$ and $v$ from $W$ we have
\begin{enumerate}
\item[(1)] $P_{u,v}(q)=0$, if $u\nleq v$;
\item[(2)] $P_{u,u}(q) =1$;
\item[(3)] $\text{deg}(P_{u,v}(q)) \leq \lfloor \frac{1}{2} (\ell(v)-\ell(u)-1) \rfloor$, if $u\lneq v$;
\item[(4)] $$q^{\ell(v)-\ell(u)} P_{u,v}(\frac{1}{q}) = \sum_{z\in [u,v]} R_{u,z}(q) P_{z,v}(q),\qquad \text{if $u\leq v$}. $$
\end{enumerate}
Also, there is a unique family of polynomials $\{P^*_{u,v}(q)\}_{u,v\in W} \subseteq \Z[q]$, called the \red{\it inverse Kazhdan-Lusztig polynomials}, satisfying the condition that 
\begin{align*}
\sum_{y\in [u,v]} (-1)^{\ell(y)-\ell(u)} P_{u,y}(q) P^*_{y,v}(q) = \delta_{u,v}\qquad \text{for all $u,v\in W,\ u\leq v$}.
\end{align*} 

In~\cite[Theorem A.2]{KazhdanLusztig1979}, the same authors showed that $X_{wB}$ is rationally smooth along the $B$-orbit $C_{vB}\subset X_{wB}$ if and only if $X_{wB}$ is rationally smooth at $e_v$ if and only if $P_{y,w}(q) =1$ holds for every $v\leq y\leq w$.
Although Kazhdan and Lusztig proved their results in characteristic $p >0$ (and using \'etale cohomology with values in the constant sheaf $\Q$), 
their results hold true in characteristic zero as well. 
Indeed, in~\cite[Theorems C and E, mixed]{Carrell1994}, Carrel (and Peterson) proved that the following statements are equivalent for an interval $[v,w]$ in $W$: 
\begin{enumerate}
\item $X_{wB}$ is rationally smooth at $e_v$, 
\item $P_{v,w}(q)=1$,
\item $P_{y,w}(q)=1$ for every $y\in [v,w]$.
\end{enumerate} 
In an unpublished manuscript, Peterson proved that if $G$ is of type $A,D$, or $E$, then a Schubert variety is rationally smooth if and only if it is smooth. 
We want to mention one more useful fact about the singularities of Richardson varieties in relation with the singularities of Schubert varieties. 
In~\cite[Corollary 2.10]{BC2012}, Billey and Coskun show that the singular locus of $X_w^v$ is given by the union 
\begin{align*}
(X_{wB}\cap X^{vB}_{sing}) \cup (X_{wB}^{sing}\cap X^{vB}),
\end{align*}
where $X^{vB}_{sing}$ and $X_{wB}^{sing}$ denote the singular locus of $X^{vB}$ and $X_{wB}$, respectively. 
%In the rest of this section, we assume that $G$ is defined over $\C$. 

We are ready to prove our Theorem~\ref{T:Smooth} from Introduction. 
Let us recall its statement for convenience. 
\medskip

Assume that $G$ is of type $A,D$, or $E$. 
Let $X_w^v$ be a toric Richardson variety with respect to the action of $T$. 
Then the interval $[v,w]$ is a Boolean lattice if and only if $X_w^v$ is a smooth toric variety.

\begin{proof}[Proof of Theorem~\ref{T:Smooth}]

($\Rightarrow$)
Towards a contradiction, let us assume that there is a singular point $x\in X_w^v$.
By the Borel fixed point theorem, there is a $T$-fixed point $e_y \in \overline{T\cdot x}\subset X_w^v$ for some $y\in [v,w]$.
This point is also singular, since the set of singular points is always closed. 
By~\cite[Corollary 2.10]{BC2012}, $y$ is either a singular point of $X^{vB}$ or $X_{wB}$. 
In the latter case, we see that $P_{y,w}(q) \neq 1$. 
But we know that~\cite[Corollary 6.8]{Brenti1994}, if $[y,w]$ is a Boolean interval, then $P_{y,w}(q) = 1$, which is absurd. 
In the former case, we see that $e_y$ is a singular point of $w_0 X_{w_0v}$.
Equivalently, $w_0 e_y$ is a singular point of $X_{w_0v}$, implying that $P_{w_0y,w_0v} (q) \neq 1$.
By the discussion after~\cite[Theorem 7.8]{Brenti1994}, we have $P_{w_0y,w_0v} (q) = P^*_{v,y}(q)$,
where $P^*_{v,y}(q)$ is the inverse Kazhdan-Lusztig polynomial. 
%\textcolor{blue}{Let $[v,y]^*$ denote the dual of interval $[v,y]$. 
%Since the interval $[v,y]$ is a Boolean lattice, it is self-dual, that is, $[v,y]\cong [v,y]^*$.(we can remove this sentence)}
Now, Brenti's~\cite[Theorem 7.8 (v)]{Brenti1994} tells us that $P^*_{v,y}(q)$ is the Stanley $g$-polynomial of $[v,y]$. 
But it is a well-known fact~\cite[Example 3.16.8]{EC1} that the $g$-polynomial of every Boolean lattice is the constant polynomial 1. 
This contradiction shows that there are no singular points in $X_w^v$. 
Hence, $X_w^v$ is smooth. 
\medskip

($\Leftarrow$) 
Let us assume that $X_w^v$ is a smooth toric Richardson variety. 
Let $\nu := \ell(w)-\ell(v)$. 
It is well-known that if $\nu =2$, then the Hasse diagram of $[v,w]$ is a diamond, hence, $[v,w]$ is a Boolean lattice. 
We proceed with the assumption that $\nu \geq 3$. 
Let $[x,y]$ be a subinterval of $[v,w]$ such that $\ell(y)-\ell(x) =3$.   
It is a well-known result~\cite{Jantzen1979} that every length 3 interval in a Weyl group is either a 2-crown or 3-crown or a 4-crown in the sense of~\cite[Figure 2.10]{BjornerBrenti}.
At the same time, since $[v,w]$ is a lattice, so is $[x,y]$.  So, $[x,y]$ cannot be a 2-crown.
We claim that $[x,y]$ cannot be a 4-crown. 
To this end, let us consider the projective embedding of $G/B$ corresponding to the line bundle $\mathcal{L}(\varpi_1+\cdots + \varpi_n)$,
where $\varpi_i$ ($1\leq i \leq n$) are the fundamental dominant weights. 
Since this embedding is smallest dimensional projective embedding of $G/B$, we refer to it as the \red{\em minimal embedding} of the flag variety. 
By restricting the minimal embedding of $G/B$ to the toric Richardson variety $X_w^v$, we obtain a smooth, very ample lattice polytope $P$ 
as in~\cite[Theorem 2.4.3]{CLS}.
In particular, we know that $P$ is a simple polytope. 
There is an anti-isomorphism between the face lattice of $P$ and the lattice of $T$-orbit closures in $X_w^v$. 
The $T$-stable curves in $X_w^v$ correspond to the 1 codimensional faces of $P$. 
Passing to the dual of $P$ we get a simplicial polytope whose vertices (resp. edges) correspond to the $T$-fixed points (resp. $T$-stable curves)
in $X_w^v$. 
Thus, the 1-skeleton of the dual polytope is precisely the Bruhat graph of the interval $[v,w]$ in the sense of~\cite{Carrell1994}.
In particular, the Bruhat graph of the interval $[x,y]$ is the 1-skeleton of a 3 dimensional simplicial polytope whose dual is a 3 dimensional simple polytope.
Let us denote this polytope by $Q_{[x,y]}$.
It has either eight vertices (if $[x,y]$ is a 3-crown) or ten vertices (if $[x,y]$ is a 4-crown).  
This means that the dual simple polytope has either eight 2-dimensional facets or ten 2-dimensional facets. 
But classification of simple polytopes in dimension 3 is well-known. 
There are only five such polytopes, called the \red{\em Platonic solids}. 
Among these five Platonic solids there is only one with eight facets, namely the octahedron.
Since there are no Platonic solids with ten facets, we see that our dual simple polytope for $[x,y]$ is the octahedron.
It follows that the Bruhat interval $[x,y]$ has 8 elements.
In other words, $[x,y]$ is a 3-crown. 
This finishes the proof of our claim. 
Now, we will finish the proof of our theorem by using a result of Grabiner~\cite[Theorem 1, Lemma 8]{Grabiner1999},
which states the following. 

Let $Q$ be a graded poset with a maximal and a minimal element.
If the following conditions are satisfied, then $Q$ is a Boolean lattice:
\begin{enumerate}
\item For every subinterval $[a,b] \subseteq Q$ of rank at least 4, the open interval $(a,b)$ is connected.
\item Every subinterval $[a,b]\subseteq Q$ of rank 3 is a Boolean lattice.
\end{enumerate}

For Bruhat-Chevalley orders, the first condition is automatically satisfied.  
In the previous paragraph, we showed that the second condition holds whenever $X_w^v$ is a smooth toric Richardson variety. 
It follows that $[v,w]$ is a Boolean lattice. 
This finishes the proof of our theorem. 
\end{proof}

\begin{Remark}
Working with the flag variety of $SL(n,\C)$, in~\cite{LMP2021}, Park, Lee, and Masuda point out that a Richardson variety $X_w^v$ is a smooth toric variety if and only if the corresponding Bruhat interval polytope is a cube. This happens precisely when the interval $[v,w]$ is a Boolean lattice. 
\end{Remark}

\section*{Acknowledgements}

Part of this research was conducted while the first author was visiting the Okinawa Institute of Science and Technology (OIST) through the Theoretical Sciences Visiting Program (TSVP).
The first author gratefully acknowledges a partial support from the Louisiana Board of Regents grant, contract no. LEQSF(2023-25)-RD-A-21. 
The second author acknowledges the National Board for Higher Mathematics (NBHM) Post Doctoral Fellowship with Ref. Number 0203/21(5)/2022-R\&D-II/10342.

\bibliographystyle{plain}
\bibliography{references}

\end{document}